\documentclass[12pt]{article}                
\usepackage{graphicx}
\usepackage{psfrag}
\usepackage{amsmath, amssymb}
\usepackage{mathptmx}
\usepackage{ntheorem}

\usepackage{setspace}
\makeatletter
\newsavebox\myboxA
\newsavebox\myboxB
\newlength\mylenA
\newcommand*\xoverline[2][0.75]{%
    \sbox{\myboxA}{$\m@th#2$}%
    \setbox\myboxB\null
    \ht\myboxB=\ht\myboxA%
    \dp\myboxB=\dp\myboxA%
    \wd\myboxB=#1\wd\myboxA
    \sbox\myboxB{$\m@th\overline{\copy\myboxB}$}
    \setlength\mylenA{\the\wd\myboxA}
    \addtolength\mylenA{-\the\wd\myboxB}%
    \ifdim\wd\myboxB<\wd\myboxA%
       \rlap{\hskip 0.5\mylenA\usebox\myboxB}{\usebox\myboxA}%
    \else
        \hskip -0.5\mylenA\rlap{\usebox\myboxA}{\hskip 0.5\mylenA\usebox\myboxB}%
    \fi}
\makeatother

\newcommand{\Av}{{\rm Av}}
\newcommand{\A}{{\cal A}}

\newcommand{\M}{{\cal M}}
\newcommand{\FF}{{\cal F}}

\newcommand{\PP}{{\cal P}}
\newcommand{\EE}{{\cal E}}

\newcommand{\LL}{{\cal L}}
\newcommand{\X}{{\cal X}}
\newcommand{\I}{{\rm I}}
\newcommand{\II}{{\cal I}}
\newcommand{\CC}{{\cal C}}

\newcommand{\Ind}{{\hspace{0.3mm}{\rm I}\hspace{0.1mm}}}
\newcommand{\eps}{{\varepsilon}}
\newcommand{\ch}{{\mbox{\rm ch}}}

\newcommand{\e}{\mathbb{E}}
\newcommand{\ep}{{\mathbb{E}_{[p]}}}
\newcommand{\epi}{{\mathbb{E}_{[p],i}}}

\newcommand{\p}{\mathbb{P}}

\newcommand{\Reals}{\mathbb{R}}
\newcommand{\Natural}{\mathbb{N}}
\newcommand{\la}{\langle}
\newcommand{\ra}{\rangle}
\newcommand\qed{\hfill\hbox{\rlap{$\sqcap$}$\sqcup$}}

\newcommand{\sbar}{{\xoverline{s}\hspace{0.1mm}}}
\newcommand{\barWp}{{\xoverline{W}\hspace{-0.6mm}}_{p}}

\newcommand{\txi}{{\tilde{\xi}}}

\newtheorem{lemma}{Lemma}
\newtheorem{theorem}{Theorem}

\theoremstyle{nonumberplain}
\newcommand\specialref{}

\begin{document}

\title{Structure of finite-RSB asymptotic Gibbs measures\\
in the diluted spin glass models}
\author{Dmitry Panchenko\thanks{Dept. of Mathematics, Texas A\&M University, panchenk@math.tamu.edu. Partially supported by NSF grant.}\\
}
\date{}
\maketitle
\begin{abstract}
We suggest a possible approach to proving the M\'ezard-Parisi formula for the free energy in the diluted spin glass models, such as diluted $K$-spin or random $K$-sat model at any positive temperature. In the main contribution of the paper, we show that a certain small modification of the Hamiltonian in any of these models forces all finite-RSB asymptotic Gibbs measures in the sense of the overlaps to satisfy the M\'ezard-Parisi ansatz for the distribution of spins. Unfortunately, what is still missing is a description of the general full-RSB asymptotic Gibbs measures. If one could show that the general case can be approximated by finite-RSB case in the right sense then one could a posteriori remove the small modification of the Hamiltonian to recover the M\'ezard-Parisi formula for the original model.
\end{abstract} 
\vspace{0.5cm}

\section{Introduction}

This paper continues to develop the line of ideas in \cite{Pspins, HEPS, AP, 1RSB}, which are all motivated by the M\'ezard-Parisi formula for the free energy in the diluted spin glass models originating in \cite{Mezard}. This formula is closely related to the original Parisi formula \cite{Parisi79, Parisi} for the free energy in the Sherrington-Kirkpatrick model \cite{SK}, but at the same time it is more complicated, because it involves a more complicated functional order parameter that encodes some very special structure of the distribution of all spins (or all multi-overlaps) rather than the distribution of the overlaps. 

An important progress was made by Franz and Leone in \cite{FL}, who showed that the M\'ezard-Parisi formula gives an upper bound on the free energy (see also \cite{PT}). The technical details of their work are very different but, clearly, inspired by the analogous earlier result of Guerra \cite{Guerra} for the Sherrington-Kirkpatrick model, which lead the to the first proof of the Parisi formula by Talagrand in \cite{TPF}. Another proof of the Parisi formula was given later in \cite{PPF}, based on the ultrametricity property for the overlaps proved in \cite{PUltra} using the Ghirlanda-Guerra identities \cite{GG} (the general idea that stability properties, such as the Aizenman-Contucci stochastic stability \cite {AC} or the Ghirlanda-Guerra identities, could imply ultrametricity is due to Arguin and Aizenman, \cite{AA}). The proof there combined the cavity method in the form of the Aizenman-Sims-Starr representation \cite{AS2} with the description of the asymptotic structure of the overlap distribution that follows from ultrametricity and the Ghirlanda-Guerra identities \cite{GG}.

The M\'ezard-Parisi ansatz in the diluted models builds upon the ultrametric Parisi ansatz in the SK model, so it is very convenient that ultrametricity for the overlaps can be obtained just as easily in the diluted models as in the SK model, simply because the Ghirlanda-Guerra identities can be proved in these models in exactly the same way, by using a small perturbation of the Hamiltonian of the mixed $p$-spin type. However, as we mention above, the M\'ezard-Parisi ansatz describes the structure of the Gibbs measure in these models in much more detail, as we shall see below.  

Some progress toward explaining the features of this ansatz beyond ultrametricity was made in \cite{HEPS, AP}, where the so-called hierarchical exchangeability of the pure states and the corresponding Aldous-Hoover representation were proved. This representation looks very similar to what one expects in the M\'ezard-Parisi ansatz, but lacks some additional symmetry. One example where this additional symmetry can be proved rigorously was given in \cite{1RSB} for the $1$-RSB asymptotic Gibbs measures in the diluted $K$-spin model, where it was obtained as a consequence of the cavity equations for spin distributions developed rigorously in \cite{Pspins}. The main contribution of this paper is to show how this result can be extended to all finite-RSB asymptotic Gibbs measures for all diluted models. Namely, we will show that one can slightly modify the Hamiltonian in such a way that all finite-RSB asymptotic Gibbs measures satisfy the M\'ezard-Parisi ansatz as a consequence of the Ghirlanda-Guerra identities and the cavity equations.

\section{Main results}

Before we can state our main results, we will need to introduce necessary notations and definitions, as well as review a number of previous results. 

Let $K\geq 1$ be an integer fixed throughout the paper. A random clause with $K$ variables will be a random function $\theta(\sigma_{1},\ldots,\sigma_{K})$ on $\{-1,+1\}^K$ symmetric in its coordinates. The main examples we have in mind are the following. 

\medskip
\noindent
{\bf Example 1.} ($K$-spin model) Given an inverse temperature parameter $\beta>0$, the random function $\theta$ is given by
$$
\theta(\sigma_1,\ldots,\sigma_K)= \beta g \sigma_1\cdots \sigma_K,
$$ 
where $g$ is a random variable, typically, standard Gaussian or Rademacher.\qed

\medskip
\noindent
{\bf Example 2.} ($K$-sat model) Given an inverse temperature parameter $\beta>0$, the random function $\theta$ is given by
$$
\theta(\sigma_1,\ldots,\sigma_K)=-\beta\prod_{j\leq K} \frac{1+J_j \sigma_j}{2},
$$
where $(J_j)_{j\geq 1}$ are i.i.d. Bernoulli random variables with $\p(J_j=\pm 1)=1/2.$\qed

\medskip

We will denote by $\theta_I$ independent copies of the function $\theta$ for various multi-indices $I$. Given a parameter $\lambda>0$, called connectivity parameter, the Hamiltonian of a diluted model is defined by
\begin{equation}
H_N^{\mathrm{model}}(\sigma) =  \sum_{k\leq \pi(\lambda N)} \theta_k(\sigma_{i_{1,k}},\ldots, \sigma_{i_{K,k}})
\label{Ham2}
\end{equation}
where $\pi(\lambda N)$ is a Poisson random variable with the mean $\lambda N$, and the coordinate indices $i_{j,k}$ are independent for different pairs $(j,k)$ and are chosen uniformly from $\{1,\ldots, N\}$. The main goal for us would be to compute the limit of the free energy
\begin{equation}
F_N^{\mathrm{model}} = \frac{1}{N}\e \log \sum_{\sigma\in \Sigma_N} \exp H_N^{\mathrm{model}}(\sigma)
\end{equation}
as $N$ goes to infinity. The formula for this limit originates in the work of M\'ezard and Parisi in \cite{Mezard}. To state how the formula looks like, we need to recall several definitions that will be used throughout the paper. 

\medskip
\noindent
\textbf{Ruelle probability cascades (RPC, \cite{Ruelle}).} Given $r\geq 1$, consider an infinitary rooted tree of depth $r$ with the vertex set
\begin{equation}
\A = \Natural^0 \cup \Natural \cup \Natural^2 \cup \ldots \cup \Natural^r,
\label{Atree}
\end{equation}
where $\Natural^0 = \{*\}$, $*$ is the root of the tree and each vertex $\alpha=(n_1,\ldots,n_p)\in \Natural^{p}$ for $p\leq r-1$ has children 
$$
\alpha n : = (n_1,\ldots,n_p,n) \in \Natural^{p+1}
$$
for all $n\in \Natural$. Each vertex $\alpha$ is connected to the root $*$ by the path
$$
* \to n_1 \to (n_1,n_2) \to\cdots\to (n_1,\ldots,n_p) = \alpha.
$$
We will denote the set of vertices in this path by 
\begin{equation}
p(\alpha) = \bigl\{*\, , n_1, (n_1,n_2),\ldots,(n_1,\ldots,n_p)  \bigr\}.
\label{pathtoleaf}
\end{equation}
We will denote by $|\alpha|$ the distance of $\alpha$ from the root, i.e. $p$ when $\alpha=(n_1,\ldots,n_p)$. We will write $\alpha \succeq \beta$ if $\beta \in p(\alpha)$ and $\alpha\succ \beta$ if, in addition, $\alpha\not =\beta$, in which case we will say that $\alpha$ is a descendant of $\beta$, and $\beta$ is an ancestor of $\alpha$. Notice that $\beta\in p(\alpha)$ if and only if $*\preceq \beta\preceq \alpha.$ The set of leaves $\Natural^r$ of $\A$ will sometimes be denoted by $\LL(\A)$. For any $\alpha, \beta\in \A$, let
\begin{equation}
\alpha\wedge\beta 
:=
 |p(\alpha) \cap p(\beta)| - 1
\label{wedge}
\end{equation}
be the number of common vertices (not counting the root $*$) in the paths from the root to the vertices $\alpha$ and  $\beta$. In other words, $\alpha \wedge \beta$ is the distance of the lowest common ancestor of $\alpha$ and $\beta$ from the root. 

Let us consider parameters
\begin{equation}
0= \zeta_{-1} <\zeta_0 <\ldots < \zeta_{r-1} <\zeta_r = 1
\label{zetas}
\end{equation}
that will appear later in the c.d.f. of the overlap in the case when it takes finitely many values (see (\ref{zetap}) below), which is the usual functional order parameter in the Parisi ansatz. For each $\alpha\in \A\setminus \Natural^r$, let $\Pi_\alpha$ be a Poisson process on $(0,\infty)$ with the mean measure $$\zeta_{p}x^{-1-\zeta_{p}}dx$$ with $p=|\alpha|$, and we assume that these processes are independent for all $\alpha$. Let us arrange all the points in $\Pi_\alpha$ in the decreasing order,
\begin{equation}
u_{\alpha 1} > u_{ \alpha 2} >\ldots >u_{\alpha n} > \ldots,
\label{us}
\end{equation}
and enumerate them using the children $(\alpha n)_{n\geq 1}$ of the vertex $\alpha$. Given a vertex $\alpha\in \A\setminus \{*\}$ and the path $p(\alpha)$ in (\ref{pathtoleaf}), we define
\begin{equation}
w_\alpha = \prod_{\beta \in p(\alpha)} u_{\beta},
\label{ws}
\end{equation}
and for the leaf vertices $\alpha \in \LL(\A) = \Natural^r$ we define
\begin{equation}
v_\alpha = \frac{w_\alpha}{\sum_{\beta\in \Natural^r} w_\beta}.
\label{vs}
\end{equation}
These are the weights of the Ruelle probability cascades. For other vertices $\alpha\in \A\setminus \LL(\A)$ we define
\begin{equation}
v_\alpha = \sum_{\beta\in \LL(\A),\,\beta\succ \alpha} v_\beta.
\label{vsall}
\end{equation}
This definition obviously implies that $v_\alpha = \sum_{n\geq 1} v_{\alpha n}$ when $|\alpha|<r$. Let us now rearrange the vertex labels so that the weights indexed by children will be decreasing. For each $\alpha\in \A\setminus \Natural^r$, let $\pi_\alpha: \Natural \to \Natural$ be a bijection such that the sequence $(v_{\alpha \pi_\alpha(n)})_{n\geq 1}$ is decreasing. Using these local rearrangements we define a global bijection $\pi: \A\to \A$ in a natural way, as follows. We let $\pi(*)=*$ and then define
\begin{equation}
\pi(\alpha n) = \pi(\alpha) \pi_{\pi(\alpha)}(n)
\label{permute}
\end{equation}
recursively from the root to the leaves of the tree. Finally, we define
\begin{equation}
V_\alpha = v_{\pi(\alpha)} \ \mbox{ for all }\  \alpha\in \A.
\label{Vs2}
\end{equation}
The weights (\ref{vs}) of the RPC will be accompanied by random fields indexed by $\Natural^r$ and generated along the tree $\A$ as follows.

\medskip
\noindent
\textbf{Hierarchical random fields.} Let $(\omega_\alpha)_{\alpha\in \A}$ be i.i.d. random variables uniform on $[0,1]$. Given a function $h: [0,1]^r \to [-1,1]$, consider a random array indexed by $\alpha= (n_1,\ldots, n_r)\in \Natural^r$,
\begin{equation}
h_\alpha = h\bigl( (\omega_\beta)_{\beta\in p(\alpha)\setminus \{*\}}\bigr)
= h\bigl(\omega_{n_1},\omega_{n_1 n_2},\ldots, \omega_{n_1\ldots n_r} \bigr).
\label{MPfop}
\end{equation}
Note that, especially, in subscripts or superscripts we will write $n_1\ldots n_r$ instead of $(n_1,\ldots, n_r)$. We will also denote by $(\omega_\alpha^I)_{\alpha\in \A}$ and
\begin{equation}
h_\alpha^I = h\bigl( (\omega^I_\beta)_{\beta\in p(\alpha)\setminus \{*\}}\bigr)
= h\bigl(\omega^I_{n_1},\omega^I_{n_1 n_2},\ldots, \omega^I_{n_1\ldots n_r} \bigr)
\label{MPfopagain}
\end{equation}
copies of the above arrays that will be independent for all different multi-indices $I.$ The function $h$ above is the second, and more complex, functional order parameter that encodes the distribution of spins inside the pure states in the M\'ezard-Parisi ansatz, as we shall see below. This way of generating the array $(h_\alpha)$ using a function $h: [0,1]^r \to [-1,1]$ is very redundant in a sense that there are many choices of the function $h$ that will produce the same array in distribution. However, if one prefers, there is a non-redundant (unique) way to encode an array of this type by a recursive tower of distributions on the set of distributions that is more common in physics literature.

\medskip
\noindent
\textbf{Extension of the definition of clause.}
Let us extend the definition of function $\theta$ on $\{-1,+1\}^K$ to $[-1,+1]^K$ as follows. Often we will need to average $\exp \theta(\sigma_1,\ldots,\sigma_K)$ over $\sigma_j \in \{-1,+1\}$ (or some subset of them) independently of each other, with some weights. If we know that the average of $\sigma_j$ is equal to $x_j\in [-1,+1]$ then the corresponding measure is given by
$$
\mu_j(\eps) = \frac{1+x_j}{2} \I(\eps=1) + \frac{1-x_j}{2} \I(\eps=-1).
$$
We would like to denote the average of $\exp \theta$ again by $\exp \theta$, which results in the definition
\begin{equation}
\exp \theta(x_1,\ldots,x_K) = \sum_{\sigma_1,\ldots,\sigma_K = \pm 1} \exp \theta(\sigma_1,\ldots,\sigma_K) \prod_{j\leq K}\mu_j(\sigma_j).
\label{Deftheta}
\end{equation}
Here is how this general definition would look like in the above two examples. In the first example of the $K$-spin model, using that $\sigma_1\cdots \sigma_K \in \{-1,+1\}$, we can write
$$
\exp \theta(\sigma_1, \ldots, \sigma_K) = \ch(\beta g)\bigl(1+\mbox{th}(\beta g)\sigma_1\cdots \sigma_K\bigr)
$$
and, clearly, after averaging,
$$
\exp \theta(x_1,\ldots, x_K) = \ch(\beta g)\bigl(1+\mbox{th}(\beta g) x_1\cdots x_K\bigr).
$$
In the second example of the $K$-sat model, using that $\prod_{j\leq K} (1+J_j \sigma_j)/2 \in\{0,1\}$, we can write
$$
\exp \theta(\sigma_1, \ldots, \sigma_K) = 
1+(e^{-\beta}-1) \prod_{j\leq K} \frac{1+J_j \sigma_j}{2}
$$
and after averaging,
$$
\exp \theta(x_1, \ldots, x_K) = 
1+(e^{-\beta}-1) \prod_{j\leq K} \frac{1+J_j x_j}{2}.
$$

\medskip
\noindent
\textbf{The M\'ezard-Parisi formula.}
Let $\pi(\lambda K)$ and $\pi(\lambda(K-1))$ be Poisson random variables with the means $\lambda K$ and $\lambda (K-1)$ correspondingly and consider
\begin{equation}
A_\alpha(\eps)
=
\sum_{k\leq \pi( \lambda K)} \theta_{k}(h_\alpha^{1,k},\ldots,h_\alpha^{K-1,k},\eps)
\label{Aibef}
\end{equation}
for $\eps\in\{-1,+1\}$ and
\begin{equation}
B_\alpha 
= 
\sum_{k\leq \pi(\lambda(K-1))} \theta_{k}(h_\alpha^{1,k},\ldots,h_\alpha^{K,k}).
\label{Bef}
\end{equation}
Let $\Av$ denote the average over $\eps=\pm 1$ and consider the following functional
\begin{equation}
\PP(r,\zeta,h)
=
\log 2 + \e \log \sum_{\alpha\in\Natural^r} v_\alpha\, \Av \exp A_\alpha(\eps)
-
\e\log \sum_{\alpha\in\Natural^r} v_\alpha \exp B_\alpha
\label{CalP}
\end{equation}
that depends on $r$, the parameters (\ref{zetas}) and the choice of the functions $h$ in (\ref{MPfop}). Then the M\'ezard-Parisi ansatz predicts that 
\begin{eqnarray}
\lim_{N\to\infty} F_N^{\mathrm{model}} = 
\inf_{r,\zeta,h} \PP(r,\zeta,h),
\label{FE}
\end{eqnarray}
at least in the above two examples of the $K$-spin and $K$-sat models. We will see below that all the parameters have a natural interpretation in terms of the structure of the Gibbs measure.

\medskip
\noindent
\textbf{Franz-Leone upper bound.}
As we mentioned in the introduction, it was proved in \cite{FL} that
\begin{eqnarray}
F_N^{\mathrm{model}} \leq \inf_{r,\zeta,h} \PP(r,\zeta,h)
\label{FL}
\end{eqnarray}
for all $N$, in the $K$-spin and $K$-sat models for even $K$. Their proof was rewritten in a slightly different language in \cite{PT} to make it technically simpler, and it was observed by Talagrand in \cite{SG} that the proof actually works for all $K\geq 1$ in the $K$-sat model. 

As a natural starting point for proving matching lower bound, a strengthened analogue of the Aizenman-Sims-Starr representation \cite{AS2} for diluted models was obtained in \cite{Pspins} in the language of the so called asymptotic Gibbs measures. We will state this representation in Theorem \ref{Th1} below for a slightly modified Hamiltonian, while also ensuring that the asymptotic Gibbs measures satisfy the Ghirlanda-Guerra identities. To state this theorem, we need to recall a few more definitions.

\medskip
\noindent
\textbf{Asymptotic Gibbs measures.} The Gibbs (probability) measure corresponding to a Hamiltonian $H_N(\sigma)$ on $\{-1,+1\}^N$ is defined by 
\begin{equation}
G_N(\sigma) = \frac{\exp H_N(\sigma)}{Z_N},
\label{GibbsN}
\end{equation} 
where the normalizing factor $Z_N=\sum_{\sigma} \exp H_N(\sigma)$ is called the partition function. To define the notion of the asymptotic Gibbs measure, we will assume that the distribution of the process $(H_N(\sigma))_{\sigma\in\{-1,+1\}^N}$ is invariant under the permutations of the coordinates of $\sigma$ - this property is called symmetry between sites, and it clearly holds in all the models we consider. 
 
Let $(\sigma^\ell)_{\ell\geq 1}$ be an i.i.d. sequence of replicas from the Gibbs measure $G_N$ and let $\mu_N$ be the joint distribution of the array $(\sigma_i^\ell)_{1\leq i\leq N, \ell\geq 1}$ of all spins for all replicas under the average product Gibbs measure $\e G_N^{\otimes \infty}$,
\begin{equation}
\mu_N\Bigl( \bigl\{\sigma_i^\ell = a_i^\ell \ :\ 1\leq i\leq N, 1\leq \ell \leq n \bigr\} \Bigr)
=
\e G_N^{\otimes n}\Bigl( \bigl\{\sigma_i^\ell = a_i^\ell \ :\ 1\leq i\leq N, 1\leq \ell \leq n \bigr\} \Bigr)
\label{muN}
\end{equation}
for any $n\geq 1$ and any $a_i^\ell \in\{-1,+1\}$.  We extend $\mu_N$ to a distribution on $\{-1,+1\}^{\Natural\times\Natural}$ simply by setting $\sigma_i^\ell=1$ for $i\geq N+1.$ Let $\M$ denote the set of all possible limits of $(\mu_N)$ over subsequences with respect to the weak convergence of measures on the compact product space $\{-1,+1\}^{\Natural\times\Natural}$. 

Because of the symmetry between sites, all measures in $\M$ inherit from $\mu_N$ the invariance under the permutation of both spin and replica indices $i$ and $\ell.$ By the Aldous-Hoover representation \cite{Aldous}, \cite{Hoover2} for such distributions, for any $\mu\in\M$, there exists a measurable function $s:[0,1]^4\to\{-1,+1\}$ such that $\mu$ is the distribution of the array 
\begin{equation}
s_i^\ell=s(w,u_\ell,v_i,x_{i,\ell}),
\label{sigma}
\end{equation}
where the random variables $w,(u_\ell), (v_i), (x_{i,\ell})$ are i.i.d. uniform on $[0,1]$. The function $s$ is defined uniquely for a given $\mu\in \M$ up to measure-preserving transformations (Theorem 2.1 in \cite{Kallenberg}), so we can identify the distribution $\mu$ of array $(s_i^\ell)$ with $s$. Since $s$ takes values in $\{-1,+1\}$, the distribution $\mu$ can be encoded by the function
\begin{equation}
\sbar(w,u,v) = \e_x\, s(w,u,v,x),
\label{fop}
\end{equation}
where $\e_x$ is the expectation in $x$ only. The last coordinate $x_{i,\ell}$ in (\ref{sigma}) is independent for all pairs $(i,\ell)$, so it plays the role of `flipping a coin' with the expected value $\sbar(w,u_\ell,v_i)$. Therefore, given the function (\ref{fop}), we can redefine $s$ by
\begin{equation}
s(w,u_\ell,v_i,x_{i,\ell}) = 2 \Ind \Bigl(x_{i,\ell} \leq \frac{1+ \sbar(w,u_\ell,v_i) }{2}\Bigr) -1
\label{sigmatos}
\end{equation}
without affecting the distribution of the array $(s_i^\ell)$. 

We can also view the function $\sbar$ in (\ref{fop}) in a more geometric way as a random  measure on the space of functions, as follows. Let $du$ and $dv$ denote the Lebesgue measure on $[0,1]$ and let us define a (random) probability measure 
\begin{equation}
G = G_w = du \circ \bigl(u\to \sbar(w,u,\cdot)\bigr)^{-1}
\label{Gibbsw}
\end{equation}
on the space of functions of $v\in [0,1]$,
\begin{equation}
H = \bigl\{ \|\sbar\|_\infty  \leq 1 \bigr\},
\label{spaceH}
\end{equation} 
equipped with the topology of $L^2([0,1], dv)$. We will denote the scalar product in $L^2([0,1], dv)$ by $h^1\cdot h^2$ and the corresponding $L^2$ norm by $\|h\|$. The random measure $G$ in (\ref{Gibbsw}) is called an {asymptotic Gibbs measure}. The whole process of generating spins can be broken into several steps:
\begin{enumerate}
\item[(i)] generate the Gibbs measure $G=G_w$ using the uniform random variable $w$; 

\item[(ii)] consider i.i.d. sequence $\sbar^{\ell} = \sbar(w,u_{\ell},\cdot)$ of replicas from $G$, which are functions in $H$; 

\item[(iii)] plug in i.i.d. uniform random variables $(v_i)_{i\geq 1}$ to obtain the array $\sbar^\ell(v_i) = \sbar(w,u_\ell,v_i)$;

\item[(iv)] finally, use this array to generate spins as in (\ref{sigmatos}). 
\end{enumerate}
For a different approach to this definition via exchangeable random measures see also \cite{Austin}. From now on, we will keep the dependence of the random measure $G$ on $w$ implicit, denote i.i.d. replicas from $G$ by $(\sbar^\ell)_{\ell\geq 1}$, which are now functions on $[0,1]$, and denote the sequence of spins (\ref{sigmatos}) corresponding to the replica $\sbar^\ell$ by
\begin{equation}
S(\sbar^\ell) = \Bigl( 2 \Ind\Bigl(x_{i,\ell} \leq \frac{1+ \sbar^\ell(v_i) }{2}\Bigr) -1 \Bigr)_{i\geq 1} \in \{-1,+1\}^\Natural.
\label{SpinsEll}
\end{equation}
Because of the geometric nature of the asymptotic Gibbs measures $G$ as measures on the subset of $L^2([0,1],dv)$, the distance and scalar product between replicas play a crucial role in the description of the structure of $G$. We will denote the scalar product between replicas $\sbar^\ell$ and $\sbar^{\ell'}$ by $R_{\ell,\ell'} = \sbar^\ell\cdot \sbar^{\ell'}$, which is more commonly called {the overlap} of $\sbar^\ell$ and $\sbar^{\ell'}$. Let us notice that the overlap $R_{\ell,\ell'}$ is a function of spin sequence (\ref{SpinsEll}) generated by $\sbar^\ell$ and $\sbar^{\ell'}$ since, by the strong law of large numbers,
\begin{equation}
R_{\ell,\ell'} = \int \! \sbar^\ell(v) \sbar^{\ell'}(v)\, dv =
\lim_{j\to\infty} \frac{1}{j}\sum_{i=1}^j S\bigl(\sbar^\ell\bigr)_i \,S\bigl(\sbar^{\ell'}\bigr)_i
\label{overlapinfty}
\end{equation}
almost surely. 

\medskip
\noindent
\textbf{The Ghirlanda-Guerra identities.} Given $n\geq 1$ and replicas $\sbar^1,\ldots, \sbar^n$, we will denote the array of spins (\ref{SpinsEll}) corresponding to these replicas by 
\begin{equation}
S^n = \bigl(S(\sbar^\ell) \bigr)_{1\leq \ell\leq n}.
\label{Sn}
\end{equation}
We will denote by $\la\,\cdot\,\ra$ the average over replicas $\sbar^\ell$ with respect to $G^{\otimes \infty}$. In the interpretation of the step (ii) above, this is the same as averaging over $(u_\ell)_{\ell\geq 1}$ in the sequence $(\sbar(w,u_{\ell},\cdot))_{\ell\geq 1}$. Let us denote by $\e$ the expectation with respect to random variables $w$, $(v_i)$ and $x_{i,\ell}$.

We will say that the measure $G$ on $H$ satisfies the Ghirlanda-Guerra identities if for any $n\geq 2,$ any bounded measurable function $f$ of the spins $S^n$ in (\ref{Sn}) and any bounded measurable function $\psi$ of one overlap,
\begin{equation}
\e \bigl\la f(S^n)\psi(R_{1,n+1}) \bigr\ra = \frac{1}{n}\hspace{0.3mm} \e\bigl\la f(S^n) \bigr\ra \hspace{0.3mm} \e\bigr\la \psi(R_{1,2})\bigr\ra + \frac{1}{n}\sum_{\ell=2}^{n}\e\bigl\la f(S^n)\psi(R_{1,\ell})\bigr\ra.
\label{GG}
\end{equation}
Another way to express the Ghirlanda-Guerra identities is to say that, conditionally on $S^n$, the law of $R_{1,n+1}$ is given by  the mixture 
\begin{equation}
\frac{1}{n} \hspace{0.3mm}\zeta + \frac{1}{n}\hspace{0.3mm} \sum_{\ell=2}^n \delta_{R_{1,\ell}},
\label{GGgen}
\end{equation}
where $\zeta$ denotes the distribution of $R_{1,2}$ under the measure $\e G^{\otimes 2}$,
\begin{equation}
\zeta(\ \cdot\ ) = \e G^{\otimes 2}\bigl(R_{1,2}\in\ \cdot\ \bigr).
\label{zeta}
\end{equation}
The identities (\ref{GG}) are usually proved for the function $f$ of the overlaps $(R_{\ell,\ell'})_{\ell,\ell'\leq n}$ instead of $S^n$, but exactly the same proof yields (\ref{GG}) as well (see e.g. Section 3.2 in \cite{SKmodel}). It is well known that these identities arise from the Gaussian integration by parts of a certain Gaussian perturbation Hamiltonian against the test function $f$, and one is free to choose this function to depend on all spins and not only overlaps.

\medskip
\noindent
\textbf{Modification of the model Hamiltonian.} Next, we will describe a crucial new ingredient that will help us classify all finite-RSB asymptotic Gibbs measures. Let us consider a sequence $(g^d)_{d\geq 1}$ of independent Gaussian random variables satisfying 
\begin{equation}
\e (g^d)^2 \leq 2^{-d}\epsilon^{\mathrm{pert}},
\label{gsvar}
\end{equation}
where $\eps^{\mathrm{pert}}$ is a fixed small parameter, and consider the following random clauses of $d$ variables,
$$
\theta^d(\sigma_1,\ldots,\sigma_d) = g^d \frac{1+\sigma_{1}}{2}\cdots \frac{1+\sigma_{d}}{2}.
$$
We will denote by $g_{I}^d$ and $\theta_{I}^d$ independent copies over different multi-indices $I$. We will define a perturbation Hamiltonian by 
\begin{equation}
H_N^{\mathrm{pert}}(\sigma)=\sum_{i\leq N} \theta_{i}^1(\sigma_i) 
+\sum_{d\geq 2} \sum_{k\leq \pi_d(N)}\theta_{k}^d(\sigma_{i_{1,d,k}},\ldots, \sigma_{i_{d,d,k}}),
\label{HNpertmain}
\end{equation}
where $\pi_d(N)$ are Poisson random variables with the mean $N$ independent over $d\geq 2$ and $i_{I}$ are chosen uniformly from $\{1,\ldots, N\}$ independently for different indices $I$. Notice that, because of (\ref{gsvar}), this Hamiltonian is well defined. We will now work with the new Hamiltonian given by
\begin{equation}
H_N(\sigma) = H_N^{\mathrm{model}}(\sigma) + H_N^{\mathrm{pert}}(\sigma).
\label{HamMain}
\end{equation}
Notice that the second term $H_N^{\mathrm{pert}}$ is of the same order as the model Hamiltonian, but its size is controlled by the parameter $\eps^{\mathrm{pert}}$. For example, if we consider the free energy
\begin{equation}
F_N = \frac{1}{N}\e \log \sum_{\sigma\in \Sigma_N} \exp H_N(\sigma)
\label{FNmod}
\end{equation}
corresponding to this modified Hamiltonian, letting $\eps^\mathrm{pert}$ go to zero will give the free energy of the original model.

Finally, as in (\ref{Deftheta}), let us extend the definition of $\theta^d$ by
\begin{equation}
\exp \theta^d(x_1,\ldots,x_d) = 1+(e^{g^d}-1) \frac{1+x_{1}}{2}\cdots \frac{1+x_{d}}{2}
\label{thetadetx}
\end{equation}
to $[-1,+1]^d$ from $\{-1,+1\}^d.$

\medskip
\noindent
\textbf{The cavity equations for the modified Hamiltonian.}
Let us now recall the cavity equations for the distribution of spins proved in \cite{Pspins}. These equations will be slightly modified here to take into account that the perturbation Hamiltonian $H_N^{\mathrm{pert}}(\sigma)$ will now also contribute to the cavity fields.

We will need to pick various sets of different spin coordinates in the array $(s_i^\ell)$ in (\ref{sigma}), and it is inconvenient to enumerate them using one index $i\geq 1$. Instead, we will use multi-indices $I= (i_1,\ldots, i_n)$ for $n\geq 1$ and $i_1,\ldots, i_n\geq 1$ and consider 
\begin{equation}
s_{I}^\ell = s(w,u_\ell, v_{I},x_{I,\ell}),
\end{equation}
where all the coordinates are uniform on $[0,1]$ and independent over different sets of indices. Similarly, we will denote
\begin{equation}
\sbar_{I}^\ell = \sbar^{\ell}(v_I) = \sbar(w,u_\ell, v_{I}).
\end{equation}
For convenience, below we will separate averaging over different replicas $\ell$, so when we average over one replica we will drop the superscript $\ell$ and simply write
\begin{equation}
\sbar_{I} = \sbar(v_I) = \sbar(w,u, v_{I}).
\label{sG}
\end{equation}
Now, take arbitrary integers $n, m, q\geq 1$ such that $n\leq m.$ The index $q$ will represent the number of replicas selected, $m$ will be the total number of spin coordinates and $n$ will be the number of cavity coordinates. For each replica index $\ell\leq q$ we  consider an arbitrary subset of coordinates $C_\ell\subseteq \{1,\ldots, m\}$ and split them into cavity and non-cavity coordinates,
\begin{equation}
C_\ell^1 = C_\ell\cap \{1,\ldots, n\},\,\,\,
C_\ell^2=C_\ell\cap \{n+1,\ldots,m\}.
\label{C12}
\end{equation}
The following quantities represent the $i$th coordinate cavity field of the modified Hamiltonian (\ref{HamMain}) in the thermodynamic limit,
\begin{align}
A_i(\eps) =& \sum_{k\leq \pi_i(\lambda K)} \theta_{k,i}( \sbar_{1,k,i}, \ldots, \sbar_{K-1,k,i}, \eps) 
\nonumber
\\
& + \theta_i^1(\eps)
+\sum_{d\geq 2} \sum_{k\leq \pi_i(d)}\theta_{k,i}^d(\sbar_{{1,d,k,i}},\ldots, \sbar_{{d-1,d,k,i}},\eps),
\label{Ai}
\end{align}
where $\pi_i(d)$ and $\pi_i(\lambda K)$ are Poisson random variables with the mean $d$ and $\lambda K$, independent of each other and independent over $d\geq 2$ and $i\geq 1$. Compared to \cite{Pspins}, now we have additional terms in the second line in (\ref{Ai}) coming from the perturbation Hamiltonian (\ref{HNpertmain}). Next, let us denote
\begin{equation}
{A}_i = \log \Av \exp {A}_i(\eps) \ \mbox{ and }\
\xi_i  = \frac{\Av \eps \exp {A}_i(\eps) }{\exp {A}_i},
\end{equation}
where $\Av$ denotes the uniform average over $\eps = \pm 1$. Recall that $\la\,\cdot\,\ra$ denotes the average with respect to the asymptotic Gibbs measure $G$. Define
\begin{equation}
{U}_\ell = \Bigl\la \prod_{i\in C_\ell^1} \xi_i \prod_{i\in C_\ell^2} \sbar_i
\,\exp \sum_{i\leq n} {A}_i \Bigr\ra \ \mbox{ and } \
{V} =\Bigl\la \exp \sum_{i\leq n} {A}_i \Bigr\ra.
\label{Ulbar2}
\end{equation}
Then we will say that an asymptotic Gibbs measure $G$ satisfies the cavity equations if
\begin{equation}
\e \prod_{\ell\leq q} \Bigl\la\, \prod_{i\in C_\ell} \sbar_i \Bigr\ra
=\e \prod_{\ell\leq q} \frac{U_\ell}{V}
\label{SC}
\end{equation}
for all choice of $n,m,q$ and sets $C_\ell^1, C_\ell^2$.

\medskip
\noindent
\textbf{The Aizenman-Sims-Starr type lower bound.}
Consider a random measure $G$ on $H$ in (\ref{spaceH}) and let $\sbar_I$ be generated by a replica $\sbar$ from this measure as in (\ref{sG}). From now on we will denote by $\pi(c)$ a Poisson random variable with the mean $c$ and we will assume that different appearances of these in the same equation are independent of each other and all other random variables. This means that if we write $\pi(a)$ and $\pi(b)$, we assume them to be independent even if $a$ happens to be equal to $b$. Consider
\begin{align}
A(\eps) =& \sum_{k\leq \pi(\lambda K)} \theta_{k}( \sbar_{1,k}, \ldots, \sbar_{K-1,k}, \eps) 
\nonumber
\\
& + \theta_i^1(\eps)
+\sum_{d\geq 2} \sum_{k\leq \pi(d)}\theta_{k}^d(\sbar_{{1,d,k}},\ldots, \sbar_{{d-1,d,k}},\eps),
\label{Aiag}
\end{align}
for $\eps\in\{-1,+1\}$ and
\begin{align}
B =& \sum_{k\leq \pi(\lambda (K-1))} \theta_{k}( \sbar_{1,k}, \ldots, \sbar_{K,k}) 
\nonumber
\\
& +\sum_{d\geq 2} \sum_{k\leq \pi(d-1)}\theta_{k}^d(\sbar_{{1,d,k}},\ldots, \sbar_{{d,d,k}}).
\label{Bag}
\end{align}
Again, compared to \cite{Pspins}, we have additional terms in the second line in (\ref{Aiag}) and (\ref{Bag}) coming from the perturbation Hamiltonian (\ref{HNpertmain}). Consider the following functional
\begin{equation}
\PP(G)
=
\log 2 + \e \log \Bigl\la \Av \exp A(\eps) \Bigr\ra
-
\e\log \Bigl\la \exp B \Bigr\ra.
\label{PPG}
\end{equation}
The following is a slight modification of the (lower bound part of the) main result in \cite{Pspins} in the setting of the diluted models.
\begin{theorem}\label{Th1}
The lower limit of the free energy in (\ref{FNmod}) satisfies
\begin{equation}
\liminf_{N\to \infty} F_N \geq \inf_G \PP(G),
\label{PPGeq}
\end{equation}
where the infimum is taken over random measures $G$ on $H$ that satisfy the Ghirlanda-Guerra identities (\ref{GG}) and the cavity equations (\ref{SC}).
\end{theorem}
We will call the measures $G$ that appear in this theorem asymptotic Gibbs measures, because that is exactly how they arise in \cite{Pspins}. The main difference from \cite{Pspins} is that we also include the requirement that the measures $G$ satisfy the Ghirlanda-Guerra identities in addition to the cavity equations. This can be ensured in exactly the same way as in the Sherrington-Kirkpatrick model by way of another small perturbation of the Hamiltonian (see e.g. \cite{HEPS}, where this was explained for the $K$-sat model). We are not going to prove Theorem \ref{Th1} in this paper, because it does not require any new ideas which are not already explained in \cite{Pspins, HEPS, 1RSB}, and the main reason we stated it here is to provide the motivation for our main result below. Of course, the proof involves some technical modifications to take into account the presence of the new perturbation term (\ref{HNpertmain}), but these are not difficult.

Instead, we will focus on the main new idea and the main new contribution of the paper, which is describing the structure of measures $G$ that satisfy the Ghirlanda-Guerra identities and the cavity equations in the case when the overlap $R_{1,2} = \sbar \cdot \sbar '$ of any two points $\sbar$ and $\sbar'$ in the support of $G$ takes finitely many, say, $r+1$ values,
\begin{equation}
0\leq q_0< q_1<\ldots< q_r \leq 1,
\label{finiteoverlap}
\end{equation}
for any $r\geq 1$ - the so called $r$-step replica symmetry breaking (or $r$-RSB) case. To state the main result, let us first recall several known consequences of the Ghirlanda-Guerra identities.

\medskip
\noindent
\textbf{Consequences of the Ghirlanda-Guerra identities.}
By Talagrand's positivity principle (see \cite{SG, SKmodel}), if the Ghirlanda-Guerra identities hold then the overlap can take only non-negative values, so the fact that the values in (\ref{finiteoverlap}) are between $0$ and $1$ is not a constraint. Another consequence of the Ghirlanda-Guerra identities (Theorem 2.15 in \cite{SKmodel}) is that with probability one the random measure $G$ is concentrated on the sphere on radius $\sqrt{q_r}$, i.e. $G(\|\sbar\|^2 = q_r)=1.$ Since we assume that the overlap takes finitely many values, $G$ is also purely atomic. Finally (see \cite{PUltra} or Theorem 2.14 in \cite{SKmodel}), with probability one the support of $G$ is ultrametric, 
$
G^{\otimes 3}(R_{2,3} \geq \min(R_{1,2},R_{1,3}))=1.
$  
By ultrametricity, for any $q_p$, the relation defined by
\begin{equation}
\sbar\sim_{q_p} \sbar' \Longleftrightarrow \sbar\cdot\sbar' \geq q_p
\label{qclusters}
\end{equation}
is an equivalence relation on the support of $G$. We will call these $\sim_q$ equivalence clusters simply $q$-clusters. Let us now enumerate all the $q_p$-clusters  defined by (\ref{qclusters}) according to Gibbs' weights as follows. Let $H_{*}$ be the entire support of $G$ so that $V_* = G(H_*) =1$. Next, the support is split into $q_1$-clusters $(H_n)_{n\geq 1}$, which are then enumerated in the decreasing order of their weights $V_n = G(H_n)$, 
\begin{equation}
V_1 > V_2 > \ldots > V_n > \ldots > 0.
\label{purelabelsfirst}
\end{equation}
We then continue recursively over $p\leq r-1$ and enumerate the $q_{p+1}$-subclusters $(H_{\alpha n})_{n\geq 1}$ of a cluster $H_\alpha$ for $\alpha\in \Natural^p$ in the non-increasing order of their weights $V_{\alpha n} = G(H_{\alpha n})$, 
\begin{equation}
V_{\alpha 1} > V_{\alpha 2} > \ldots > V_{\alpha n} > \ldots > 0.
\label{purelabels}
\end{equation}
Thus, all these clusters were enumerated $(V_\alpha)_{\alpha\in \A}$ by the vertices of the tree $\A$ in (\ref{Atree}). It is not a coincidence that we used the same notation as in (\ref{Vs2}). It is another well-known consequence of the Ghirlanda-Guerra identities that the distribution of these weights coincides with the reordering of the weights of the Ruelle probability cascades as in (\ref{Vs2}) with the parameters (\ref{zetas}) given by
\begin{equation}
\zeta_p = \e G^{\otimes 2}\bigl( R_{1,2} \leq q_p \bigr)
\label{zetap}
\end{equation}
for $p=0,\ldots,r.$ The $q_r$-clusters are the points of the support of $G$ - these are called pure states. They were enumerated by $\alpha\in\Natural^r$ and, if we denote them by $\sbar_\alpha$, \begin{equation} 
G(\sbar_\alpha) = V_\alpha \,\,\mbox{ for }\,\, \alpha\in \Natural^r.
\label{Gdiscrete}
\end{equation}
Recall that we generate the array $\sbar_i^\ell$ (or $\sbar_I^\ell$ for general index $I$) by first sampling replicas $\sbar^\ell$ from the measure $G$ (which are functions on $[0,1]$) and then plugging in i.i.d. uniform random variables $v_i$, i.e. $\sbar_i^\ell = \sbar^\ell(v_i)$. In the discrete setting (\ref{Gdiscrete}), this is equivalent to sampling $\alpha$ according to the weights $V_\alpha$ and then plugging in $v_i$ into $\sbar_\alpha,$ i.e.  $\sbar_\alpha(v_i)$. Therefore, in order to describe the distribution of the array $(\sbar^\ell(v_i))_{i,\ell\geq 1}$ it is sufficient to describe the joint distribution of the arrays
$$
\bigl(V_{\alpha}\bigr)_{\alpha\in\Natural^r} \,\,\mbox{ and }\,\,
\bigl(\sbar_{\alpha}(v_i)\bigr)_{i\geq 1, \alpha\in\Natural^r}.
$$
In addition to the fact that $(V_\alpha)$ corresponds to some reordering of weights of the Ruelle probability cascades, it was proved in \cite{AP, HEPS} that if the measure $G$ satisfies the Ghirlanda-Guerra identities then (see Theorem $1$ and equation (36) and (37) in \cite{HEPS}):
\begin{enumerate}
\item[(i)] the arrays $(V_{\alpha})_{\alpha\in\Natural^r}$ and $(\sbar_{\alpha}(v_i))_{i\geq 1, \alpha\in\Natural^r}$ are independent;

\item[(ii)] there exists a function $h:[0,1]^{2(r+1)}\to[-1,1]$ such that 
\begin{equation}
\Bigl(\sbar_\alpha(v_i) \Bigr)_{i\geq 1, \alpha\in\Natural^r}
\, \stackrel{d}{=}\,
\Bigl( h\bigl( (\omega_\beta)_{\beta\in p(\alpha)}, (\omega_\beta^i)_{\beta\in p(\alpha)} \bigr) \Bigr)_{i\geq 1,\alpha\in\Natural^r},
\label{sigmaf2}
\end{equation}
where, as above, $\omega_\alpha$ and $\omega_\alpha^i$ for $\alpha\in\A$ are i.i.d. uniform random variables on $[0,1]$. 
\end{enumerate}
The M\'ezard-Parisi ansatz predicts that in the equation (\ref{sigmaf2}) one can replace the function $h$ by a function that does not depend on the coordinates $(\omega_\beta)_{\beta\in p(\alpha)}$, which would produce exactly the same fields as in (\ref{MPfopagain}). We will show  that this essentially holds for finite-RSB asymptotic Gibbs measures.

\medskip
\noindent
\textbf{Consequence of the cavity equations.} The main result of the paper is the following.

\begin{theorem}\label{Th2}
If a random measure $G$ on $H$ in (\ref{spaceH}) satisfies the Ghirlanda-Guerra identities (\ref{GG}) and the cavity equations (\ref{SC}) and the overlap takes $r+1$ values in (\ref{finiteoverlap}) then there exists a function $h:[0,1]^{r+2}\to[-1,1]$ such that 
\begin{equation}
\Bigl(\sbar_\alpha(v_i) \Bigr)_{i\geq 1, \alpha\in\Natural^r}
\, \stackrel{d}{=}\,
\Bigl( h\bigl( \omega_*, (\omega_\beta^i)_{\beta\in p(\alpha)} \bigr) \Bigr)_{i\geq 1,\alpha\in\Natural^r},
\label{sigmaf3}
\end{equation}
where $\omega_*$ and $\omega_\alpha^i$ for $\alpha\in\A$ are i.i.d. uniform random variables on $[0,1]$. 
\end{theorem}
In other words, the cavity equations (\ref{SC}) allow us to simplify (\ref{sigmaf2}) and to get rid of the dependence on the coordinates $\omega_\beta$ for $\beta\in \A\setminus \{*\}$. Notice that, compared to the M\'ezard-Parisi ansatz, we still have the dependence on $\omega_*$ in (\ref{sigmaf3}). However, from the point of view of computing the free energy this is not an issue at all, because the average in $\omega_*$ is on the outside of the logarithm in (\ref{PPG}) and when we minimize over $G$ in (\ref{PPGeq}), we can replace the average over $\omega_*$ by the infimum. Of course, the infimum over $G$ in (\ref{PPGeq}) could involve measures that are not of finite-RSB type, and this is the main obstacle to finish the proof of the M\'ezard-Parisi formula, if this approach can be made to work. 

If one could replace the infimum in (\ref{PPGeq}) over measures $G$ that satisfy the finite-RSB condition in (\ref{finiteoverlap}) (in addition to the cavity equations and the Ghirlanda-Guerra identities) then, using Theorem \ref{Th2} and replacing the average over $\omega_*$ by the infimum, we get the lower bound that essentially matches the Franz-Leone upper bound, except that now we have additional terms in the second line in (\ref{Aiag}) and (\ref{Bag}) compared to (\ref{Aibef}) and (\ref{Bef}) coming from the perturbation Hamiltonian (\ref{HNpertmain}). However, these terms are controlled by $\eps^\mathrm{pert}$ in (\ref{gsvar}) and, letting it go to zero, one could remove the dependence of the lower bound on these terms and match the Franz-Leone upper bound.

\section{General idea of the proof}\label{Sec2ilabel}

The main goal of this paper is to show that the function $h$ that generates the array $\sbar_i^\alpha$ in (\ref{sigmaf2}), 
$$
\sbar_i^\alpha = h\bigl( (\omega_\beta)_{\beta\in p(\alpha)}, (\omega_\beta^i)_{\beta\in p(\alpha)} \bigr), 
$$
can be replaced by a function that does not depend on the coordinates $\omega_\beta$ for $*\prec \beta\preceq \alpha$. We will show this by induction, removing one coordinate at a time from the leaf $\alpha$ up to the root $*$. Our induction assumption will be the following: for $p\in \{0,\ldots, r-1\}$, suppose that, instead of (\ref{sigmaf2}), the array $\sbar_i^\alpha$ for $i\geq 1,\alpha\in \Natural^r$, is generated by
\begin{equation}
\sbar_i^\alpha = h\bigl( (\omega_\beta)_{\beta\in p(\alpha), |\beta|\leq p+1}, (\omega_\beta^i)_{\beta\in p(\alpha)} \bigr)
\label{indass}
\end{equation}
for some function $h$ that does not depend on the coordinates $\omega_\beta$ for $|\beta|\geq p+2$. Notice that this holds for $p=r-1$, and we would like to show that one can replace $h$ by $h'$ that also does not depend on $\omega_{\beta}$ for $|\beta|=p+1$, without affecting the distribution of the array $\sbar_i^\alpha$. 

Often, we will work with a subtree of $\A$ that `grows out' of a vertex at the distance $p$ or $p+1$ from the root, which means that all paths from the root to the vertices in that subtree pass through this vertex. In that case, for certainty, we will fix the vertex to be $[p] = (1,2,\ldots,p)$ or $[p+1]$. We will denote by $\ep$ the expectation with respect to the random variables $\omega_\beta$, $\omega_\beta^i$ indexed by the descendants of $[p]$, i.e. $[p]\prec \beta$, and by $\epi$ the expectation with respect to $\omega_\beta^i$ for $[p]\prec \beta$. Our goal will be to prove the following.
\begin{theorem}\label{Sec6iTh1}
Under the assumption (\ref{indass}), for any $k\geq 1$ and any $[p+1]\preceq \alpha_1,\ldots,\alpha_k \in \Natural^r$, the expectation $\epi \sbar_i^{\alpha_1}\cdots \sbar_i^{\alpha_k}$ with respect to $(\omega_\beta^i)_{[p+1]\preceq \beta}$ does not depend on $\omega_{[p+1]}$ almost surely.
\end{theorem}
Here $i\geq 1$ is arbitrary but fixed and $\alpha_1,\ldots,\alpha_k$ need not be different, so the quantities $\epi \sbar_i^{\alpha_1}\cdots \sbar_i^{\alpha_k}$ represent all possible joint moments with respect to $(\omega_\beta^i)_{[p+1]\preceq \beta}$ of the random variables $\sbar_i^\alpha$ for $[p+1]\preceq \alpha\in \Natural^r$.
It will take us the rest of the paper to prove this result, and right now we will only explain why it completes the induction step. 

The reason is identical to the situation of an exchangeable sequence $X_n = h(\omega,\omega_n)$ (say, bounded in absolute value by one) such that all moments $\e_{\omega_1} X_1^k$ for $k\geq 1$ with respect to $\omega_1$ do not depend on $\omega$. In this case if we choose any function $h'(\omega_1)$ with this common set of moments then the sequences $(h'(\omega_n))$ and $(X_n)$ have the same distribution, which can be seen by comparing their joint moments. For example, we can choose $h'(\,\cdot\,) = h(\omega^*,\,\cdot\,)$ for any $\omega^*$ from the set of measure one on which all moments $\e_{\omega_1} X_1^k$ coincide with their average values $\e X_1^k$. 

We can do the same in the setting of Theorem \ref{Sec6iTh1}, which can be rephrased as follows: for almost all $(\omega_\beta,\omega_\beta^i)_{\beta\preceq [p]}$ and $\omega_{[p+1]}$,
$$
\epi \sbar_i^{\alpha_1}\cdots \sbar_i^{\alpha_k} = \ep \sbar_i^{\alpha_1}\cdots \sbar_i^{\alpha_k}
$$
for all $k\geq 1$ and $[p+1]\preceq \alpha_1,\ldots,\alpha_k \in \Natural^r$, where $\ep$ now also includes the average in $\omega_{[p+1]}.$ This means that we can find $\omega_{[p+1]}=\omega_{[p+1]}^*$
such that the equality of all these moments holds for almost all $(\omega_\beta,\omega_\beta^i)_{\beta\preceq [p]}$. If we now set
$$
h'\bigl( (\omega_\beta)_{\beta\in p(\alpha), |\beta|\leq p}, (\omega_\beta^i)_{\beta\in p(\alpha)} \bigr)
=
h\bigl( (\omega_\beta)_{\beta\in p(\alpha), |\beta|\leq p}, \omega_{[p+1]}^*,  (\omega_\beta^i)_{\beta\in p(\alpha)} \bigr)
$$
then by comparing the joint moments one can see that
\begin{equation}
\bigl(\sbar_i^\alpha\bigr)_{i\geq 1,\alpha\in\Natural^r} \stackrel{d}{=} 
\Bigl(h'\bigl( (\omega_\beta)_{\beta\in p(\alpha), |\beta|\leq p}, (\omega_\beta^i)_{\beta\in p(\alpha)} \bigr) \Bigr)_{i\geq 1,\alpha\in\Natural^r},
\end{equation}
which completes the induction step. 

The proof of Theorem \ref{Sec6iTh1} will proceed by a certain induction on the shape of the configuration $\alpha_1,\ldots,\alpha_k$, where by the shape of the configuration we essentially mean the matrix $(\alpha_\ell\wedge \alpha_{\ell'})_{\ell,\ell'\leq k}$ (or its representation by a tree that consists of all paths $p(\alpha_\ell)$). It is clear that the quantity
$$
\epi \sbar_i^{\alpha_1}\cdots \sbar_i^{\alpha_k}
$$
depends on $[p+1]\preceq \alpha_1,\ldots,\alpha_k \in \Natural^r$ only through the shape. The induction will be somewhat involved and we will explain exactly how it will work toward the end of the paper, once we have all the tools ready. However, we need to mention now that the induction will have an important \textit{monotonicity property}: whenever we have proved the statement of Theorem \ref{Sec6iTh1} for some $\alpha_1,\ldots,\alpha_k$, we have also proved it for any subset of these vertices. At this moment, we will suppose that $[p+1]\preceq \alpha_1,\ldots,\alpha_k \in \Natural^r$ are such that the following holds:
\begin{enumerate}
\item[(M)] For any subset $S\subseteq \{1,\ldots, k\}$, $\epi \prod_{\ell\in S}\sbar_i^{\alpha_\ell}$ does not depend on $\omega_{[p+1]}$.
\end{enumerate} 
Then over the next sections we will obtain some implications of this assumption using the cavity equations. Finally, in the last section we will show how to use these implications inductively to prove Theorem \ref{Sec6iTh1} for any choice of $\alpha_1,\ldots,\alpha_k$. Of course, the starting point of the induction will be the case of $k=1$ that we will obtain first. In fact, in this case the statement will be even stronger and will not assume that (\ref{indass}) holds (i.e. we only assume (\ref{sigmaf2})).
\begin{lemma}\label{Sec2iLem1}
For any $p=0,\ldots, r-1$ and any $[p]\prec \alpha \in \Natural^r$, the expectation $\epi \sbar_i^{\alpha}$ with respect to $(\omega_\beta^i)_{[p]\prec \beta}$ does not depend on $(\omega_\beta)_{[p]\prec \beta}$ almost surely.
\end{lemma}
\textbf{Proof.} Consider $[p]\prec \alpha, \beta\in\Natural^r$ such that $\alpha\wedge \beta = p$. By (\ref{sigmaf2}), it is clear that the overlap of two pure states satisfies
\begin{equation}
R_{\alpha, \beta}:= \sbar_\alpha\cdot \sbar_{\beta} 
=
\int_0^1\! \sbar_\alpha(v) \sbar_{\beta}(v) \, dv
\stackrel{d}{=}
\e_i \, \sbar_i^{\alpha} \sbar_i^{\beta},
\label{Sec2eq1}
\end{equation}
where $\e_i$ denotes the expectation in random variables $\omega_\eta^i$ that depend on the spin index $i$. By construction, we enumerated the pure states $\sbar_\alpha$ in (\ref{Gdiscrete}) so that
\begin{equation}
R_{\alpha,\beta} = q_{\alpha\wedge\beta}
\label{RabSec2}
\end{equation}
and, since $\alpha\wedge \beta = p$, we get that, almost surely,
$$
q_p = R_{\alpha,\beta} = \e_i \sbar_i^{\alpha}  \sbar_i^{\beta} 
= \e_i \bigl(\epi \sbar_i^{\alpha} \epi \sbar_i^{\beta}\bigr).
$$
If we denote $v = (\omega_\eta)_{\eta\preceq [p]}$, $v_1 = (\omega_\eta)_{[p]\prec \eta\preceq \alpha}$, $v_2 = (\omega_\eta)_{[p]\prec \eta\preceq \beta}$ and $u = (\omega_\eta^i)_{\eta\preceq [p]}$ then
$$
\epi \sbar_i^{\alpha} =\varphi(v,v_1,u),\,\, 
\epi \sbar_i^{\beta}=\varphi(v,v_2,u)
$$
for some function $\varphi,$ the random variables $v,v_1,v_2,u$ are independent, and the above equation can be written as
$$
q_p = \e_u \varphi(v,v_1,u) \varphi(v,v_2,u)
$$
for almost all $v,v_1,v_2$. This means that for almost all $v$, the above equality holds for almost all $v_1, v_2$. Let us fix any such $v$ and let $\mu_v$ be the image of the Lebesgue measure on $[0,1]^{r-p}$ by the map $v_1 \to \varphi(v,v_1,\cdot)\in L^2([0,1]^{p+1},du)$. Then the above equation means that, if we sample independently two points from $\mu_v$, with probability one their scalar product in $L^2$ will be equal to $q_p$. This can happen only if the measure $\mu_v$ is concentrated on one point in $L^2$, which means that the function $\varphi$ does not depend on $v_1$.
\qed

\medskip
\noindent
Before we start using the cavity equations, we will explain a property of the Ruelle probability cascades that will play the role of the main technical tool throughout the paper.

\section{Key properties of the Ruelle probability cascades}\label{Sec2label}

The property described in this section will be used in two ways - directly, in order to obtain some consequences of the cavity equations, and indirectly, as a representation tool to make certain computations possible. This property is proved in Theorem 4.4 in \cite{SKmodel} in a more general form, but here we will need only a special case as follows.

Let us consider a random variable $z$ taking values in some measurable space $(\X, {\cal B})$ (in our case, this will always be some nice space, such as $[0,1]^n$ with the Borel $\sigma$-algebra) and let $z_\alpha$ be its independent copies indexed by the vertices of the tree $\alpha\in \A\setminus \{*\}$ excluding the root. Recall the parameters $(\zeta_p)_{0\leq p\leq r-1}$ in (\ref{zetas}). Let us consider a measurable bounded function 
\begin{equation}
X_r:\X^r \to \Reals
\label{Xr}
\end{equation}
and, recursively over $0\leq p\leq r-1$, define functions $X_p:\X^p \to\Reals$ by
\begin{equation}
X_p(x) = \frac{1}{\zeta_p}\log \e_z \exp \zeta_p X_{p+1}(x,z),
\label{ch54Xp}
\end{equation}
where the expected value $\e_z$ is with respect to $z$. In particular, $X_0$ is a constant. Let us define
\begin{equation}
W_p(x,y) = \exp \zeta_p\bigl(X_{p+1}(x,y) - X_p(x)\bigr)
\label{ch31Wp}
\end{equation}
for $x\in \X^p$ and $y\in\X.$ Let us point out that, by the definition (\ref{ch54Xp}), $\e_z W_p(x,z) = 1$ and, therefore, for each $x\in\X^p,$ we can think of $W_p(x,\, \cdot\, )$ as a change of density that yields the following conditional distribution on $\X$ given $x\in\X^p$,
\begin{equation}
\nu_p(x,B) = \e_z W_p(x,z)\, \I(z\in B).
\label{ch51transitionprime}
\end{equation}
For $p=0$, $\nu_0$ is just a probability distribution on $(\X,{\cal B}).$ Let us now generate the array $\tilde{z}_\alpha$ for $\alpha\in \A\setminus \{*\}$ iteratively from the root to the leaves as follows. Let $\tilde{z}_n$ for $n\in\Natural$ be i.i.d. random variables with the distribution $\nu_0$. If we already constructed $\tilde{z}_\alpha$ for $|\alpha|\leq p$ then, given any $\alpha\in \Natural^p$, we generate $\tilde{z}_{\alpha n}$ independently for $n\geq 1$ from the conditional distribution $\nu_p((z_\beta)_{*\prec \beta\preceq \alpha},\, \cdot\,)$, and these are generated independently over different such $\alpha.$ Notice that the distribution of the array $(\tilde{z}_\alpha)_{\alpha\in\A\setminus\{*\}}$ depends on the distribution of $z$, function $X_r$ and parameters $(\zeta_p)_{0\leq p\leq r-1}$.

With this definition, the expectation ${\e} f((\tilde{z}_\alpha)_{\alpha\in F})$ for a finite subset $F\subset \A\setminus\{*\}$ can be written as follow. For $\alpha\in\A\setminus \{*\}$, let
\begin{equation}
W_\alpha = W_{|\alpha|-1}\bigl( (z_\beta)_{*\prec \beta\prec \alpha}, z_\alpha \bigr).
\label{Sec2Walpha}
\end{equation}
Slightly abusing notation, we could also write this simply as $W_\alpha = W_{|\alpha|-1}( (z_\beta)_{*\prec \beta\preceq \alpha})$. Given a finite subset ${C} \subset \A\setminus\{*\}$, let
$$
p({C}) = \bigcup_{\alpha\in {C}} p(\alpha) \setminus \{*\}
\,\,\mbox{ and }\,\,
W_{{C}} = \prod_{\alpha\in p({C})} W_\alpha.
$$
Then, the above definition of the array $(\tilde{z}_\alpha)$ means that
\begin{equation}
{\e} f\bigl((\tilde{z}_\alpha)_{\alpha\in {C}} \bigr)
=
\e W_{{C}}  f\bigl(({z}_\alpha)_{\alpha\in {C}}\bigr).
\label{Sec2expectW}
\end{equation}
Simply, to average over $(\tilde{z}_\alpha)_{\alpha\in {C}}$ we need to use changes of density over all vertices in the paths from the root leading to the vertices $\alpha\in {C}$.

The meaning of the above construction will be explained by the following result. Recall the Ruelle probability $(v_\alpha)_{\alpha\in \Natural^r}$ cascades in (\ref{vs}) and define new random weights on $\Natural^r$,
\begin{equation}
\tilde{v}_\alpha
=
\frac{v_\alpha \exp X_r((z_\beta)_{*\prec \beta\preceq \alpha})}{\sum_{\alpha\in \Natural^r}v_\alpha \exp X_r((z_\beta)_{*\prec \beta\preceq \alpha})},
\label{tildevs}
\end{equation}
by the change of density proportional to $\exp X_r((z_\beta)_{*\prec \beta\preceq \alpha})$. We will say that a bijection $\pi:\A\to\A$ of the vertices of the tree $\A$ preserves the parent-child relationship if children $\alpha n$ of $\alpha$ are mapped into children of $\pi(\alpha)$, $\pi(\alpha n) = (\pi(\alpha),k)$ for some $k\in\Natural$. Another way to write this is to say that $\pi(\alpha)\wedge \pi(\beta) = \alpha\wedge\beta$ for all $\alpha,\beta\in \A.$ For example, the bijection $\pi$ defined in (\ref{permute}), (\ref{Vs2}), is of this type. Theorem 4.4 in \cite{SKmodel} gives the following generalization of the Bolthausen-Sznitman invariance property for the Poisson-Dirichlet point process (Proposition A.2 in \cite{Bolthausen}). 
\begin{theorem}
There exists a random bijection $\rho:\A\to\A$ of the vertices of the tree $\A$, which preserves the parent-child relationship, such that
\begin{equation}
\bigl(\tilde{v}_{\rho(\alpha)}\bigr)_{\alpha\in\Natural^r}
\stackrel{d}{=}
\bigl(v_\alpha\bigr)_{\alpha\in\Natural^r},\,\,
\bigl(z_{\rho(\alpha)} \bigr)_{\alpha\in \A\setminus \{*\}}
\stackrel{d}{=}
\bigl(\tilde{z}_\alpha \bigr)_{\alpha\in \A\setminus \{*\}}
\label{Th3eq}
\end{equation}
and these two arrays are independent of each other. 
\end{theorem}
This result will be more useful to us in a slightly different formulation in terms of the sequence $(V_\alpha)_{\alpha\in \Natural^r}$ in (\ref{Vs2}). Namely, if we denote by
\begin{equation}
\tilde{V}_\alpha
=
\frac{V_\alpha \exp X_r((z_\beta)_{*\prec \beta\preceq \alpha})}{\sum_{\alpha\in \Natural^r}V_\alpha \exp X_r((z_\beta)_{*\prec \beta\preceq \alpha})}
\label{tildeVs}
\end{equation}
then the following holds.
\begin{theorem}\label{Th4label}
There exists a random bijection $\rho:\A\to\A$ of the vertices of the tree $\A$, which preserves the parent-child relationship, such that
\begin{equation}
\bigl(\tilde{V}_{\rho(\alpha)}\bigr)_{\alpha\in\Natural^r}
\stackrel{d}{=}
\bigl(V_\alpha\bigr)_{\alpha\in\Natural^r},\,\,
\bigl(z_{\rho(\alpha)} \bigr)_{\alpha\in \A\setminus \{*\}}
\stackrel{d}{=}
\bigl(\tilde{z}_\alpha \bigr)_{\alpha\in \A\setminus \{*\}}
\label{Th4eq}
\end{equation}
and these two arrays are independent of each other. 
\end{theorem}
\textbf{Proof.} We have to apply twice the following simple observation. Suppose that we have a random array $(v_\alpha')_{\alpha\in\Natural^r}$ of positive weights that add up to one and array $(z_\alpha')_{\alpha\in\A\setminus\{*\}}$ generated along the tree similarly to $(\tilde{z}_\alpha)$ above - namely, ${z}_n'$ for $n\in\Natural$ are i.i.d. random variables with some distribution $\nu_0$ and, if we already constructed ${z}_\alpha'$ for $|\alpha|\leq p$ then, given any $\alpha\in \Natural^p$, we generate ${z}_{\alpha n}'$ independently for $n\geq 1$ from some conditional distribution $\nu_p((z_\beta)_{*\prec \beta\preceq \alpha},\, \cdot\,)$, and these are generated independently over different such $\alpha.$ Suppose that $(v_\alpha')_{\alpha\in\Natural^r}$ and $(z_\alpha')_{\alpha\in\A\setminus\{*\}}$ are independent. Consider any random permutation $\rho:\A\to\A$ that preserves the parent-child relationship, which depends only on $(v_\alpha')_{\alpha\in\Natural^r}$, i.e. it is a measurable function of this array. Then the arrays 
$$
\bigl(v_{\rho(\alpha)}' \bigr)_{\alpha\in\Natural^r} \,\,\mbox{ and }\,\,  \bigl(z_{\rho(\alpha)}' \bigr)_{\alpha\in\A\setminus\{*\}}
$$ 
are independent and
$$
\bigl(z_{\rho(\alpha)}' \bigr)_{\alpha\in\A\setminus\{*\}}
\stackrel{d}{=}
\bigl(z_\alpha' \bigr)_{\alpha\in\A\setminus\{*\}}.
$$
This is obvious because, conditionally on $\rho$, the array $z_{\rho(\alpha)}'$ is generated exactly like $z_\alpha'$ along the tree, so its conditional distribution does not depend on $\rho$. One example of such permutation $\rho$ is the permutation defined in (\ref{vsall}), (\ref{permute}), (\ref{Vs2}), that sorts the cluster weights indexed by $\alpha\in \A\setminus \Natural^r$ defined by
\begin{equation}
v_\alpha' = \sum_{\alpha\prec \beta\in\Natural^r} v_\beta'.
\label{vsallprime}
\end{equation}
Namely, for each $\alpha\in \A\setminus \Natural^r$, we let $\pi_\alpha: \Natural \to \Natural$ be a bijection such that the sequence $v_{\alpha \pi_\alpha(n)}'$ is decreasing for $n\geq 1$ (we assume that all these cluster weights are different as is the case for the Ruelle probability cascades), let $\pi(*)=*$ and define
\begin{equation}
\pi(\alpha n) = \pi(\alpha) \pi_{\pi(\alpha)}(n)
\label{permuteprime}
\end{equation}
recursively from the root to the leaves of the tree. Let us denote
\begin{equation}
\mathrm{Sort}\Bigl(
\bigl(v_\alpha' \bigr)_{\alpha\in\Natural^r}, \bigl(z_\alpha' \bigr)_{\alpha\in\A\setminus\{*\}}
\Bigr)
:=
\Bigl(
\bigl(v_{\pi(\alpha)}' \bigr)_{\alpha\in\Natural^r}, \bigl(z_{\pi(\alpha)}' \bigr)_{\alpha\in\A\setminus\{*\}}
\Bigr).
\label{sort}
\end{equation}
Notice that this sorting operation depends only on $(v_\alpha')_{\alpha\in\Natural^r}$, so it does not affect the distribution of $(z_\alpha')_{\alpha\in\A\setminus\{*\}}$. Now, let us show how (\ref{Th3eq}) implies (\ref{Th4eq}). First of all, the permutation $\rho$ in the equation (\ref{Th4eq}) is just the sorting operation described above,
$$
\Bigl(
\bigl(\tilde{V}_{\rho(\alpha)}\bigr)_{\alpha\in\Natural^r},
\bigl(z_{\rho(\alpha)} \bigr)_{\alpha\in \A\setminus \{*\}}
\Bigr)
=
\mathrm{Sort}
\Bigl(
\bigl(\tilde{V}_{\alpha}\bigr)_{\alpha\in\Natural^r},
\bigl(z_{\alpha} \bigr)_{\alpha\in \A\setminus \{*\}}
\Bigr).
$$
Let $\pi$ be the permutation in (\ref{permute}), (\ref{Vs2}) and, trivially,
$$
\mathrm{Sort}
\Bigl(
\bigl(\tilde{V}_{\alpha}\bigr)_{\alpha\in\Natural^r},
\bigl(z_{\alpha} \bigr)_{\alpha\in \A\setminus \{*\}}
\Bigr)
=
\mathrm{Sort}
\Bigl(
\bigl(\tilde{V}_{\pi^{-1}(\alpha)}\bigr)_{\alpha\in\Natural^r},
\bigl(z_{\pi^{-1}(\alpha)} \bigr)_{\alpha\in \A\setminus \{*\}}
\Bigr),
$$
since the sorting operation does not depend on how we index the array. On the other hand, by the definition (\ref{tildeVs}) and the fact that $V_{\pi^{-1}(\alpha)} = v_\alpha$,
$$
\tilde{V}_{\pi^{-1}(\alpha)}
=
\frac{v_\alpha \exp X_r((z_{\pi^{-1}(\beta)})_{\beta\in p(\alpha)})}{\sum_{\alpha\in \Natural^r}v_\alpha \exp X_r((z_{\pi^{-1}(\beta)})_{\beta\in p(\alpha)})}.
$$
Also, since the permutation $\pi$ depends only on $(v_\alpha)$, by the above observation, the arrays $(v_{\alpha})_{\alpha\in\Natural^r}$ and  $(z_{\pi^{-1}(\alpha)})_{\alpha\in\A\setminus\{*\}}$ are independent and
$$
\bigl(z_{\pi^{-1}(\alpha)} \bigr)_{\alpha\in\A\setminus\{*\}}
\stackrel{d}{=}
\bigl(z_\alpha \bigr)_{\alpha\in\A\setminus\{*\}}.
$$
Comparing with the definition (\ref{tildevs}), this gives that
$$
\Bigl(
\bigl(\tilde{V}_{\pi^{-1}(\alpha)}\bigr)_{\alpha\in\Natural^r},
\bigl(z_{\pi^{-1}(\alpha)} \bigr)_{\alpha\in \A\setminus \{*\}}
\Bigr)
\stackrel{d}{=} 
\Bigl(
\bigl(\tilde{v}_{\alpha}\bigr)_{\alpha\in\Natural^r},
\bigl(z_{\alpha} \bigr)_{\alpha\in \A\setminus \{*\}}
\Bigr),
$$
and all together we showed that
$$
\Bigl(
\bigl(\tilde{V}_{\rho(\alpha)}\bigr)_{\alpha\in\Natural^r},
\bigl(z_{\rho(\alpha)} \bigr)_{\alpha\in \A\setminus \{*\}}
\Bigr)
\stackrel{d}{=} 
\mathrm{Sort}
\Bigl(
\bigl(\tilde{v}_{\alpha}\bigr)_{\alpha\in\Natural^r},
\bigl(z_{\alpha} \bigr)_{\alpha\in \A\setminus \{*\}}
\Bigr).
$$
Since we already use the notation $\rho$, let us denote the permutation $\rho$ in (\ref{Th3eq}) by $\rho'$. Then (\ref{Th3eq}) implies
\begin{align*}
&
\mathrm{Sort}
\Bigl(
\bigl(\tilde{v}_{\alpha}\bigr)_{\alpha\in\Natural^r},
\bigl(z_{\alpha} \bigr)_{\alpha\in \A\setminus \{*\}}
\Bigr)
=
\mathrm{Sort}
\Bigl(
\bigl(\tilde{v}_{\rho'(\alpha)}\bigr)_{\alpha\in\Natural^r},
\bigl(z_{\rho'(\alpha)} \bigr)_{\alpha\in \A\setminus \{*\}}
\Bigr)
\\
&
\stackrel{d}{=} 
\mathrm{Sort}
\Bigl(
\bigl({v}_{\alpha}\bigr)_{\alpha\in\Natural^r},
\bigl(\tilde{z}_{\alpha} \bigr)_{\alpha\in \A\setminus \{*\}}
\Bigr)
=
\Bigl(
\bigl({V}_{\alpha}\bigr)_{\alpha\in\Natural^r},
\bigl(\tilde{z}_{\pi(\alpha)} \bigr)_{\alpha\in \A\setminus \{*\}}
\Bigr).
\end{align*}
Finally, since the sorting permutation $\pi$ depends only on the array $(v_\alpha)$, by the above observation, the array $(\tilde{z}_{\pi(\alpha)})$ is independent of $(V_\alpha)$ and has the same distribution as $(\tilde{z}_{\alpha})$. This finishes the proof.
\qed

\section{Cavity equations for the pure states}\label{Sec3label}

In this section, we will obtain some general consequences of the cavity equations (\ref{SC}) that do not depend on any inductive assumptions. In the next section, we will push this further under the assumption (M) made in Section \ref{Sec2ilabel}. First of all, let us rewrite the cavity equations (\ref{SC}) taking into account the consequences of the Ghirlanda-Guerra identities in (\ref{Gdiscrete}) and (\ref{sigmaf2}). Let us define
\begin{align}
\sbar^\alpha_I \ =\ &  \ h\bigl( (\omega_\beta)_{\beta\in p(\alpha)}, (\omega_\beta^I)_{\beta\in p(\alpha)} \bigr),
\label{sialpha}
\\
A_i^\alpha(\eps) \ =\ & \sum_{k\leq \pi_i(\lambda K)} \theta_{k,i}( \sbar^\alpha_{1,k,i}, \ldots, \sbar^\alpha_{K-1,k,i}, \eps) 
\nonumber
\\
&  \ + \theta_i^1(\eps)
+\sum_{d\geq 2} \sum_{k\leq \pi_i(d)}\theta_{k,i}^d(\sbar^\alpha_{{1,d,k,i}},\ldots, \sbar^\alpha_{{d-1,d,k,i}},\eps),
\label{Aieps}
\\
A_i^\alpha \ =\ & \log \Av \exp A_i^\alpha(\eps),
\label{Aialpha}
\\
\xi_i^\alpha \ =\ & \frac{\Av \eps \exp A_i^\alpha(\eps) }{\exp A_i^\alpha},
\label{Sec3xiialpha}
\end{align}
and let $A^\alpha = \sum_{i\leq n} A_i^\alpha$. We will keep the dependence of $A^\alpha$ on $n$ implicit for simplicity of notation. Then (\ref{Ulbar2}) can be redefined by (using equality in distribution (\ref{sigmaf2}))
\begin{equation}
{U}_\ell =\, \sum_{\alpha\in \Natural^r} V_\alpha \prod_{i\in C_\ell^1} \xi_i^\alpha \prod_{i\in C_\ell^2} \sbar_i^\alpha\,\exp A^\alpha \ \mbox{ and } \
{V} =\, \sum_{\alpha\in \Natural^r } V_\alpha \exp A^\alpha.
\label{Ulbaralpha}
\end{equation}
Moreover, if we denote
\begin{equation}
\tilde{V}_\alpha = \frac{V_\alpha \exp A^\alpha}{{V}}
= \frac{V_\alpha \exp A^\alpha}{\sum_{\alpha\in \Natural^r} V_\alpha \exp A^\alpha}
\label{Valpha}
\end{equation}
then the cavity equations (\ref{SC}) take form
\begin{equation}
\e \prod_{\ell\leq q} \sum_{\alpha\in \Natural^r} V_\alpha \prod_{i\in C_\ell} \sbar_i^\alpha
=\e \prod_{\ell\leq q} \sum_{\alpha\in \Natural^r} \tilde{V}_\alpha \prod_{i\in C_\ell^1} \xi_i^\alpha \prod_{i\in C_\ell^2} \sbar_i^\alpha.
\label{SCbaralpha}
\end{equation}
We can also write this  as
\begin{equation}
\e \sum_{\alpha_1,\ldots, \alpha_q} V_{\alpha_1}\cdots V_{\alpha_q} \prod_{\ell\leq q}
\prod_{i\in C_\ell} \sbar_i^{\alpha_\ell}
=\e \sum_{\alpha_1,\ldots, \alpha_q} \tilde{V}_{\alpha_1} \cdots \tilde{V}_{\alpha_q} \prod_{\ell\leq q} \prod_{i\in C_\ell^1} \xi_i^{\alpha_\ell} \prod_{i\in C_\ell^2} \sbar_i^{\alpha_\ell}.
\label{SCnew}
\end{equation}
We will now use this form of the cavity equations to obtain a different form directly for the pure states that does not involve averaging over the pure states. 

Let us formulate the main result of this section. Let $\FF$ be the $\sigma$-algebra generated by the random variables that are not indexed by $\alpha\in \A\setminus \{*\}$, namely, 
\begin{equation}
\theta_{k,i}, \theta_i^1, \theta^d_{k,i}, \pi_i(\lambda K), \pi_i(d), \omega_*, \omega_*^I
\label{FF}
\end{equation}
for various indices, excluding the random variables $\omega_{\alpha}$ and $\omega_{\alpha}^I$ that are indexed by $\alpha\in \A\setminus \{*\}$. Let $\II_i$ be the set of indices $I$ that appear in various $\sbar_I^\alpha$ in (\ref{Aieps}), i.e. $I$ of the type $(\ell, k,i,)$ or $(\ell, d,k,i)$. Let $\II = \cup_{i\geq 1} \II_i$ and let
\begin{equation}
z_\alpha^i = \bigl(\omega_\alpha^I \bigr)_{I\in\II_i},\,\,
z_\alpha =  \bigl( \omega_\alpha, (z_\alpha^i)_{i\geq 1} \bigr)= \bigl( \omega_\alpha, (\omega_\alpha^I)_{I\in\II} \bigr),
\label{zeealpha}
\end{equation}
Notice that with this notation, conditionally on $\FF$, the random variables $A_i^\alpha$ and $\xi_i^\alpha$ in (\ref{Aialpha}), (\ref{Sec3xiialpha}) for $\alpha\in\Natural^r$ can be written as
\begin{equation}
\xi_i^\alpha = \xi_i \bigl( (\omega_\beta,z_\beta^i)_{*\prec \beta \preceq \alpha}\bigr),\,\, 
A_i^\alpha = \chi_i \bigl( (\omega_\beta, z_\beta^i)_{*\prec \beta \preceq \alpha}\bigr). 
\label{xichi}
\end{equation}
for some function $\xi_i$ and $\chi_i$ (that implicitly depend on the random variables in (\ref{FF})) and
$$
A^\alpha = X\bigl( (z_\beta)_{*\prec \beta \preceq \alpha}\bigr) := \sum_{i\leq n} \chi_i \bigl( (\omega_\beta, z_\beta^i)_{*\prec \beta \preceq \alpha}\bigr).
$$ 
In the setting of the previous section, let $X_r = X$ in (\ref{Xr}) and let $(\tilde{z}_\alpha)_{\alpha\in\A\setminus \{*\}}$ be the array generated along the tree using the conditional probabilities (\ref{ch51transitionprime}). Recall that this means the following. The definition in (\ref{ch54Xp}) can be written as
\begin{equation}
X_{|\alpha|-1}\bigl( (z_\beta)_{*\prec \beta \prec \alpha}\bigr)
=\frac{1}{\zeta_{|\alpha|-1}} \log \e_{z_\alpha} \exp \zeta_{|\alpha|-1} X_{|\alpha|}\bigl( (z_\beta)_{*\prec \beta \preceq \alpha}\bigr),
\label{Sec3Xp}
\end{equation}
where $\e_{z_\alpha}$ is the expectation in $z_\alpha$, and the definition in (\ref{ch31Wp}) can be written as
\begin{equation}
W_{|\alpha|-1} \bigl( (z_\beta)_{*\prec \beta \prec \alpha},z_\alpha\bigr) 
= 
\exp \zeta_{|\alpha|-1}\Bigl(X_{|\alpha|}\bigl( (z_\beta)_{*\prec \beta \preceq \alpha}\bigr)  - X_{|\alpha|-1}\bigl( (z_\beta)_{*\prec \beta \prec \alpha}\bigr) \Bigr).
\label{Sec3Walpha}
\end{equation}
Then the array $(\tilde{z}_\alpha)_{\alpha\in\A\setminus \{*\}}$ is generated along the tree from the root to the leaves according to the conditional probabilities in (\ref{ch51transitionprime}), namely, given $(\tilde{z}_\beta)_{*\prec \beta \prec \alpha}$ we generated $\tilde{z}_\alpha$ by the change of density $W_{|\alpha|-1} \bigl( (\tilde{z}_\beta)_{*\prec \beta \prec \alpha},\ \cdot \ \bigr)$. Let us emphasize one more time that this entire construction is done conditionally on $\FF$. Also, notice that the coordinates $\omega_\alpha^I$ in (\ref{zeealpha}) were independent for different $I$, but the corresponding coordinates $\tilde{\omega}_\alpha^I$ of
$
\tilde{z}_\alpha = \bigl( \tilde{\omega}_\alpha, (\tilde{\omega}_\alpha^I)_{I\in\II} \bigr)
$
are no longer independent, because $X_r$ and the changes of density $W_p$ depend on all of them. As in (\ref{zeealpha}) and (\ref{xichi}), let us denote
\begin{equation}
\tilde{z}_\alpha^i = \bigl(\tilde{\omega}_\alpha^I \bigr)_{I\in\II_i},\,\,
\tilde{z}_\alpha =  \bigl( \tilde{\omega}_\alpha, (\tilde{z}_\alpha^i)_{i\geq 1} \bigr),\,\,
\tilde{\xi}_i^\alpha = \xi_i \bigl( (\tilde{\omega}_\beta, \tilde{z}_\beta^i)_{*\prec \beta \preceq \alpha}\bigr).
\label{Sec4tildas}
\end{equation}
We will prove the following.
\begin{theorem}\label{Sec4Th}
The equality in distribution holds (not conditionally on $\FF$),
\begin{equation}
\bigl(\txi_i^{\alpha} \bigr)_{i\leq n,\alpha\in \Natural^r}
\stackrel{d}{=}
\bigl(\sbar_i^{\alpha} \bigr)_{i\leq n,\alpha\in \Natural^r}.
\label{Sec4ThEq}
\end{equation}
\end{theorem}
\textbf{Proof.} 
As in (\ref{Sec2eq1}) and (\ref{RabSec2}), we can write
$$
R_{\alpha, \beta}= \sbar_\alpha\cdot \sbar_{\beta} 
=
\int_0^1\! \sbar_\alpha(v) \sbar_{\beta}(v) \, dv
=
\e_i \, \sbar_i^{\alpha} \sbar_i^{\beta},
$$
where $\e_i$ denotes the expectation in random variables $\omega_\beta^i$ in (\ref{sialpha}) that depend on the spin index $i$, and $R_{\alpha,\beta} = q_{\alpha\wedge\beta}$.

In the cavity equations (\ref{SCnew}), let us now make a special choice of the sets $C_\ell^2$. For each pair $(\ell,\ell')$ of replica indices such that $1\leq \ell<\ell'\leq q$, take any integer $n_{\ell,\ell'}\geq 0$ and consider a set $C_{\ell,\ell'}\subseteq \{n+1,\ldots,m\}$ of cardinality $|C_{\ell,\ell'}|=n_{\ell,\ell'}$. Let all these sets be disjoint, which can be achieved by taking $m=n+\sum_{1\leq \ell<\ell'\leq q} n_{\ell,\ell'}.$ For each $\ell\leq q$, let
$$
C_\ell^2 = \Bigl(\bigcup_{\ell'>\ell} C_{\ell,\ell'}\Bigr) \bigcup \Bigl(\bigcup_{\ell'<\ell} C_{\ell',\ell}\Bigr).
$$
Then a given spin index $i\in \{n+1,\ldots,m\}$ appears in exactly two sets, say, $C_\ell^2$ and $C_{\ell'}^2$, and the expectation of (\ref{SCnew}) in $(\omega_{\beta}^i)$ will produce a factor $\e_i \,s_i^{\alpha_\ell} s_i^{\alpha_{\ell'}} = R_{\alpha_\ell,\alpha_{\ell'}}$. For each pair $(\ell,\ell')$, there will be exactly $n_{\ell,\ell'}$ such factors, so averaging in (\ref{SCnew}) in the random variables $(\omega_{\beta}^i)$ for all $i\in \{n+1,\ldots,m\}$ will result in
\begin{equation}
\e \sum_{\alpha_1,\ldots, \alpha_q} V_{\alpha_1}\cdots V_{\alpha_q} \prod_{\ell<\ell'} R_{\alpha_\ell, \alpha_{\ell'}}^{n_{\ell,\ell'}}
\prod_{\ell\leq q} \prod_{i\in C_\ell^1} \sbar_i^{\alpha_\ell}
=\e \sum_{\alpha_1,\ldots, \alpha_q} \tilde{V}_{\alpha_1} \cdots \tilde{V}_{\alpha_q}  \prod_{\ell<\ell'} R_{\alpha_\ell, \alpha_{\ell'}}^{n_{\ell,\ell'}} \prod_{\ell\leq q} \prod_{i\in C_\ell^1} \xi_i^{\alpha_\ell}.
\label{SCagain}
\end{equation}
Approximating by polynomials, we can replace $\prod_{\ell<\ell'} R_{\alpha_\ell, \alpha_{\ell'}}^{n_{\ell,\ell'}}$ by the indicator of the set
\begin{equation}
\CC = \bigl\{(\alpha_1,\ldots, \alpha_q) \ | \ R_{\alpha_\ell, \alpha_{\ell'}} = q_{\ell,\ell'} \mbox{ for all } 1\leq \ell<\ell' \leq q\bigr\}
\end{equation}
for any choice of constraints $q_{\ell,\ell'}$ taking values in $\{q_0,\ldots,q_r\}$. Therefore, (\ref{SCagain}) implies
\begin{equation}
\sum_{(\alpha_1,\ldots, \alpha_q)\in \CC} \e V_{\alpha_1}\cdots V_{\alpha_q} \prod_{\ell\leq q} \prod_{i\in C_\ell^1} \sbar_i^{\alpha_\ell}
= \sum_{(\alpha_1,\ldots, \alpha_q)\in \CC} \e \tilde{V}_{\alpha_1} \cdots \tilde{V}_{\alpha_q}  \prod_{\ell\leq q} \prod_{i\in C_\ell^1} \xi_i^{\alpha_\ell}.
\label{SCF}
\end{equation}
Using the property (i) above the equation (\ref{sigmaf2}), which as we mentioned is the consequence of the Ghirlanda-Guerra identities, we can rewrite the left hand side as
$$
\sum_{(\alpha_1,\ldots, \alpha_q)\in \CC} \e V_{\alpha_1}\cdots V_{\alpha_q} \,
\e \prod_{\ell\leq q} \prod_{i\in C_\ell^1} \sbar_i^{\alpha_\ell}.
$$
Moreover, it is obvious from the definition of the array $\sbar_i^\alpha$ in (\ref{sialpha}) that the second expectation depends on $(\alpha_1,\ldots, \alpha_q)\in \CC$ only through the overlap constraints $(q_{\ell,\ell'})$, or $(\alpha_\ell\wedge \alpha_{\ell'})$. 

On the other hand, on the right hand side of (\ref{SCF}) both $\tilde{V}_\alpha$ and $\xi_i^\alpha$ depend on the same random variables through the function $A_i^\alpha(\eps)$. If we compare (\ref{tildeVs}) and (\ref{Valpha}) and apply Theorem \ref{Th4label} conditionally on $\FF$, we see that that there exists a random bijection $\rho:\A\to\A$ of the vertices of the tree $\A$ which preserves the parent-child relationship and such that
\begin{equation}
\bigl(\tilde{V}_{\rho(\alpha)}\bigr)_{\alpha\in\Natural^r}
\stackrel{d}{=}
\bigl(V_\alpha\bigr)_{\alpha\in\Natural^r},\,\,
\bigl( (\xi_i^{\rho(\alpha)} )_{i\leq n} \bigr)_{\alpha\in \A\setminus \{*\}}
\stackrel{d}{=}
\bigl( (\txi_i^{\alpha})_{i\leq n} \bigr)_{\alpha\in \A\setminus \{*\}}
\end{equation}
and these two arrays are independent of each other (all these statement are conditionally on $\FF$). If we denote by $\e'$ the conditional expectation given $\FF$ then this implies that
\begin{align*}
&
\sum_{(\alpha_1,\ldots, \alpha_q)\in \CC} \e'\, \tilde{V}_{\alpha_1} \cdots \tilde{V}_{\alpha_q}  \prod_{\ell\leq q} \prod_{i\in C_\ell^1} \xi_i^{\alpha_\ell}
=
\sum_{(\alpha_1,\ldots, \alpha_q)\in \CC} \e'\, \tilde{V}_{\rho(\alpha_1)} \cdots \tilde{V}_{\rho(\alpha_q)}  \prod_{\ell\leq q} \prod_{i\in C_\ell^1} \xi_i^{\rho(\alpha_\ell)}
\\
&=
\sum_{(\alpha_1,\ldots, \alpha_q)\in \CC} \e'\, V_{\alpha_1} \cdots V_{\alpha_q}  \prod_{\ell\leq q} \prod_{i\in C_\ell^1} \txi_i^{\alpha_\ell}
=
\sum_{(\alpha_1,\ldots, \alpha_q)\in \CC} \e\,' V_{\alpha_1} \cdots V_{\alpha_q} \,
\e'\, \prod_{\ell\leq q} \prod_{i\in C_\ell^1} \txi_i^{\alpha_\ell}.
\end{align*}
Since the distribution of $(V_\alpha)_{\alpha\in\Natural^r}$ does not depend on the condition and $\e\,' V_{\alpha_1} \cdots V_{\alpha_q} = \e V_{\alpha_1} \cdots V_{\alpha_q}$, taking the expectation gives
$$
\sum_{(\alpha_1,\ldots, \alpha_q)\in \CC} \e \tilde{V}_{\alpha_1} \cdots \tilde{V}_{\alpha_q}  \prod_{\ell\leq q} \prod_{i\in C_\ell^1} \xi_i^{\alpha_\ell}
=
\sum_{(\alpha_1,\ldots, \alpha_q)\in \CC} \e V_{\alpha_1} \cdots V_{\alpha_q} \,
\e \prod_{\ell\leq q} \prod_{i\in C_\ell^1} \txi_i^{\alpha_\ell}.
$$
This proves that
$$
\sum_{(\alpha_1,\ldots, \alpha_q)\in \CC} \e V_{\alpha_1}\cdots V_{\alpha_q} \,
\e \prod_{\ell\leq q} \prod_{i\in C_\ell^1} \sbar_i^{\alpha_\ell}
=
\sum_{(\alpha_1,\ldots, \alpha_q)\in \CC} \e V_{\alpha_1} \cdots V_{\alpha_q} \,
\e \prod_{\ell\leq q} \prod_{i\in C_\ell^1} \txi_i^{\alpha_\ell}.
$$
Again, the second expectation in the sum on the right depends on $(\alpha_1,\ldots, \alpha_q)\in \CC$ only through the overlap constraints $(q_{\ell,\ell'})$ and, since the choice of the constraints was arbitrary, we get
\begin{equation}
\e \prod_{\ell\leq q} \prod_{i\in C_\ell^1} \sbar_i^{\alpha_\ell}
=
\e \prod_{\ell\leq q} \prod_{i\in C_\ell^1} \txi_i^{\alpha_\ell}
\end{equation}
for any $\alpha_1,\ldots, \alpha_q \in \Natural^r$. Clearly, one can express any joint moment of the elements in these two arrays by choosing $q\geq 1$ large enough and choosing $\alpha_1,\ldots, \alpha_q$ and the sets $C_\ell^1$ properly, so the proof is complete.
\qed

\section{A consequence of the cavity equations for the pure states}\label{Sec4label}

We will continue using the notation of the previous section, only in this section we will take $n=2$ in Theorem \ref{Sec4Th}. Let us recall the assumption (M) made at the end of Section \ref{Sec2ilabel}: we consider some $[p+1]\preceq \alpha_1,\ldots,\alpha_k \in \Natural^r$ such that the following holds:
\begin{enumerate}
\item[(M)] For any subset $S\subseteq \{1,\ldots, k\}$, $\epi \prod_{\ell\in S}\sbar_i^{\alpha_\ell}$ does not depend on $\omega_{[p+1]}$.
\end{enumerate} 
In this section, we will obtain a further consequence of the cavity equations using that $\epi \prod_{\ell\leq k}\sbar_i^{\alpha_\ell}$ does not depend on $\omega_{[p+1]}$, but a similar consequence will hold for any subset of these vertices.

Let us denote by $\ep$ the expectation with respect to the random variables $\omega_\eta, \omega_\eta^I$ indexed by the ancestors $\eta\succ [p]$ of $[p]$. We will use the same notation $\ep$ to denote the expectation with respect to the random variables $\tilde{\omega}_\eta, \tilde{\omega}_\eta^I$ for $\eta\succ [p]$ conditionally on $\tilde{\omega}_\eta, \tilde{\omega}_\eta^I$ for $\eta\preceq [p]$ and all other random variables that generate the $\sigma$-algebra $\FF$ in (\ref{FF}). Given any finite set ${C}\subset \Natural^r$, let us denote
$$
\sbar_i^{C} = \prod_{\alpha\in {C}} \sbar^\alpha_i,\,\,
\txi_i^{C} = \prod_{\alpha\in {C}} \txi^\alpha_i
\,\,\mbox{ and }\,\,
\xi_i^{C} = \prod_{\alpha\in {C}} \xi^\alpha_i.
$$
Then the following holds for $[p+1]\preceq \alpha_1,\ldots,\alpha_k \in \Natural^r$.
\begin{lemma}\label{Sec5Lem1}
If ${C} = \{\alpha_1,\ldots,\alpha_k\}$ and $\epi \sbar_i^{C}$ does not depend on $\omega_{[p+1]}$ then
\begin{equation}
\ep  \txi_1^{{C}}  \txi_2^{{C}} = \ep  \txi_1^{{C}} \ep  \txi_2^{{C}}
\label{Sec3Lem1eq}
\end{equation}
almost surely.
\end{lemma}
\textbf{Proof.} First of all,
$$
\ep \sbar_1^{{C}} \sbar_2^{{C}} =
\ep \prod_{i\leq 2 } \epi \sbar_i^{{C}}
= \prod_{i\leq 2 } \epi \sbar_i^{{C}}
$$
almost surely, since $\epi \sbar_i^{C}$ does not depend on $\omega_{[p+1]}$. Similarly,
$$
\ep \sbar_i^{{C}} =
\ep \epi \sbar_i^{{C}}
= \epi \sbar_i^{{C}}
$$
almost surely and, therefore,
\begin{align}
0
= &\
\e \bigl(\ep \sbar_1^{{C}} \sbar_2^{{C}} - \ep \sbar_1^{{C}} \ep \sbar_2^{{C}} \bigr)^2
\nonumber
\\
= &\
\e \bigl(\ep \sbar_1^{{C}} \sbar_2^{{C}}\bigr)^2
-2 \e \bigl(\ep \sbar_1^{{C}} \sbar_2^{{C}}\bigr) \bigl(\ep \sbar_1^{{C}} \ep \sbar_2^{{C}}\bigr)
+\e \bigl(\ep \sbar_1^{{C}} \ep \sbar_2^{{C}}\bigr)^2.
\label{Sec3above}
\end{align}
Let us now rewrite each of these terms using replicas. Let ${C}_1 = {C}$ and for $j=2,3,4$ let ${C}_j = \{\alpha^j_1,\ldots,\alpha^j_k\}$ for arbitrary $[p+j]\preceq \alpha^j_1,\ldots,\alpha^j_k \in \Natural^r$ such that $\alpha^j_\ell \wedge \alpha^j_{\ell'} = \alpha_\ell \wedge \alpha_{\ell'}$ for any $\ell,\ell'\leq k.$ In other words, ${C}_j$ are copies of ${C}$ that consists of the descendants of different children of $[p]$. Therefore, we can write 
$$
\ep \sbar_1^{{C}_j} \sbar_2^{{C}_j}
=
\ep \sbar_1^{{C}_{j'}} \sbar_2^{{C}_{j'}}
\,\,\mbox{ and }\,\,
\ep \sbar_i^{{C}_j} 
=
\ep \sbar_i^{{C}_{j'}}
$$
almost surely for any $j,j'\leq 4$ and
\begin{align*}
\e \bigl(\ep \sbar_1^{{C}} \sbar_2^{{C}}\bigr)^2
= &\
\e \sbar_1^{{C}_1} \sbar_2^{{C}_1}\sbar_1^{{C}_2} \sbar_2^{{C}_2},
\\
\e \bigl(\ep \sbar_1^{{C}} \sbar_2^{{C}}\bigr) \bigl(\ep \sbar_1^{{C}} \ep \sbar_2^{{C}}\bigr)
= &\
\e \sbar_1^{{C}_1} \sbar_2^{{C}_1}\sbar_1^{{C}_2} \sbar_2^{{C}_3},
\\
\e \bigl(\ep \sbar_1^{{C}} \ep \sbar_2^{{C}}\bigr)^2
= &\
\e \sbar_1^{{C}_1} \sbar_2^{{C}_2}\sbar_1^{{C}_3} \sbar_2^{{C}_4}.
\end{align*}
By Theorem \ref{Sec4Th}, this and (\ref{Sec3above}) imply that
$$
\e \txi_1^{{C}_1} \txi_2^{{C}_1}\txi_1^{{C}_2} \txi_2^{{C}_2}
-2 \e \txi_1^{{C}_1} \txi_2^{{C}_1}\txi_1^{{C}_2} \txi_2^{{C}_3}
+\e \txi_1^{{C}_1} \txi_2^{{C}_2}\txi_1^{{C}_3} \txi_2^{{C}_4} = 0.
$$
Repeating the above computation backwards for $\txi$ instead of $\sbar$ gives
$$
\e \bigl(\ep  \txi_1^{{C}}  \txi_2^{{C}} - \ep  \txi_1^{{C}} \ep  \txi_2^{{C}} \bigr)^2 = 0
$$
and this finishes the proof.
\qed

\bigskip
\noindent
By analogy with (\ref{Sec2Walpha}) and (\ref{Sec2expectW}), let us rewrite the expectation $\ep$ with respect to the random variables $\tilde{\omega}_\alpha, \tilde{\omega}_\alpha^I$ for $\alpha\succ [p]$ in terms of the expectation with respect to the random variables ${\omega}_\alpha, {\omega}_\alpha^I$ for $\alpha\succ [p]$, writing explicitly the changes of density
\begin{equation}
W_\alpha = W_{|\alpha|-1} \bigl( (z_\beta)_{*\prec \beta \prec \alpha},z_\alpha\bigr). 
\label{Sec5Wa}
\end{equation}
As in Lemma \ref{Sec5Lem1}, let ${C} = \{\alpha_1,\ldots,\alpha_k\}$ for some $[p+1]\preceq \alpha_1,\ldots,\alpha_k \in \Natural^r$, let
$$
p([p],{C}) = \bigl\{\beta \ \bigr|\ [p+1]\preceq \beta\preceq \alpha, \alpha\in{C} \bigr\}
$$
and define
\begin{equation}
W_{[p],{C}} = \prod_{\alpha\in p([p],{C})} W_\alpha.
\label{Sec5WpC}
\end{equation}
With this notation, we can rewrite (\ref{Sec3Lem1eq}) as
\begin{equation}
\ep  \xi_1^{{C}}  \xi_2^{{C}} W_{[p],{C}} = \ep  \xi_1^{{C}}W_{[p],{C}}\, \ep  \xi_2^{{C}}W_{[p],{C}}
\label{Sec5eq1}
\end{equation}
almost surely. Notice that in (\ref{Sec3Lem1eq}) almost surely meant for almost all random variables in (\ref{FF}) that generate the $\sigma$-algebra $\FF$ and for almost all $\tilde{z}_\alpha$ for $\alpha\preceq [p]$ that are generated conditionally on $\FF$ according to the changes of density in (\ref{Sec3Walpha}). However, even though in (\ref{Sec5eq1}) we simply expressed the expectation with respect to $\tilde{z}_\alpha$ for $[p]\prec \alpha$ using the changes of density explicitly, after this averaging both sides depend on ${z}_\alpha$ for $\alpha\preceq [p]$, so almost surely now means for almost all random variables in (\ref{FF}) that generate the $\sigma$-algebra $\FF$ and for almost all ${z}_\alpha$ for $\alpha\preceq [p]$. The reason we can do this is very simple. Notice that $A_i^\alpha(\eps)$ in (\ref{Aieps}) can be bounded by
$$
 |A_i^\alpha(\eps)| \leq  c_i:= \sum_{k\leq \pi_i(\lambda K)} \|\theta_{k,i}\|_{\infty} + |g_i^1| +\sum_{d\geq 2} \sum_{k\leq \pi_i(d)}|g_{k,i}^d|,
$$
which, by the assumption (\ref{gsvar}), is almost surely finite (notice also that $c_i$ are $\FF$-measurable). By induction in (\ref{Sec3Xp}), all $|X_{|\alpha|}|\leq c=c_1+c_2$ almost surely and, therefore, all changes of density in (\ref{Sec3Walpha}) satisfy $e^{-2c}\leq W_{|\alpha|}\leq e^{2c}$ almost surely. Therefore, conditionally on $\FF$, the distribution of all $z_\alpha$ and $\tilde{z}_\alpha$ are absolutely continuous with respect to each other and, therefore, we can write almost surely equality in (\ref{Sec5eq1}) in terms of the random variables $z_\alpha$ for $\alpha\preceq [p]$.

Next, we will reformulate (\ref{Sec5eq1}) using the assumption (\ref{indass}). To simplify the notation, let us denote for any $\alpha\in \A\setminus\{*\}$,
\begin{equation}
\omega_{\preceq\alpha} =  (\omega_\beta)_{*\prec \beta\preceq \alpha},\,\,
z_{\preceq\alpha}^i = (z_\beta^i)_{*\prec \beta \preceq \alpha},\,\,
z_{\preceq\alpha} = (z_\beta)_{*\prec \beta \preceq \alpha}
\end{equation}
and define $\omega_{\prec\alpha}, z_{\prec\alpha}^i $ and $z_{\prec\alpha}$ similarly. Then, we can rewrite (\ref{xichi}) for $[p+1]\preceq \alpha\in \Natural^r$ as
\begin{equation}
\xi_i^{\alpha} = \xi_i \bigl( \omega_{\preceq [p+1]}, z^i_{\preceq \alpha}\bigr),\,\,
A_i^{\alpha} =\chi_i \bigl( \omega_{\preceq [p+1]}, z^i_{\preceq \alpha}\bigr). 
\label{Sec5xichi}
\end{equation}
Since in the previous section we set $n=2$, we have
$$
A^{\alpha} = \sum_{i\leq 2} \chi_i \bigl( \omega_{\preceq [p+1]}, z^i_{\preceq \alpha}\bigr).
$$ 
Because of the absence of the random  variables $\omega_\alpha$ for $[p+1]\prec \alpha$, the integration in $z_\alpha$ in the recursive definition (\ref{Sec3Xp}) will decouple when $[p+1]\prec \alpha$ into integration over $z_\alpha^1$ and $z_\alpha^2$. Namely, let $\chi_{i,r}=\chi_i$ and, for $[p+1] \preceq \alpha$, let us define by decreasing induction on $|\alpha|$,
\begin{equation}
\chi_{i,|\alpha|-1}\bigl( \omega_{\preceq [p+1]}, z^i_{\prec \alpha }\bigr)
=
\frac{1}{\zeta_{|\alpha|-1}} \log \e_{z_{\alpha}^i} \exp \zeta_{|\alpha|-1} \chi_{i,|\alpha|}\bigl( \omega_{\preceq [p+1]}, z^i_{\preceq \alpha}\bigr).
\label{chij}
\end{equation}
First of all, for $[p+1]\prec \alpha$, by decreasing induction on $|\alpha|$,
\begin{align}
X_{|\alpha|-1}\bigl( z_{\prec \alpha}\bigr)
= &\
\frac{1}{\zeta_{|\alpha|-1}} \log \e_{z_{\alpha}} \exp \zeta_{|\alpha|-1} X_{|\alpha|}\bigl( z_{\preceq \alpha}\bigr)
\nonumber
\\
\{\mbox{induction assumption}\}= &\
\frac{1}{\zeta_{|\alpha|-1}} \log \e_{z_{\alpha}} \exp \zeta_{|\alpha|-1}
\sum_{i\leq 2}
\chi_{i,|\alpha|}\bigl( \omega_{\preceq [p+1]}, z^i_{\preceq \alpha}\bigr)
\nonumber
\\
\{\mbox{independence}\}= &\
\frac{1}{\zeta_{|\alpha|-1}} \log \prod_{i\leq 2} \e_{z_{\alpha}^i} \exp \zeta_{|\alpha|-1}
\chi_{i,|\alpha|}\bigl( \omega_{\preceq [p+1]}, z^i_{\preceq \alpha}\bigr)
\nonumber
\\
\{\mbox{definition (\ref{chij})}\}= &\
\sum_{i\leq 2}
\chi_{i,|\alpha|-1}\bigl( \omega_{\preceq [p+1]}, z^i_{\prec \alpha}\bigr).
\label{Xichii}
\end{align}
When we do the same computation for $\alpha = [p+1]$, the expectation $\e_{z_{[p+1]}}$ also involves $\omega_{[p+1]}$, so we end up with
\begin{align}
X_{p}\bigl( z_{\preceq [p]}\bigr)
= &\
\frac{1}{\zeta_p} \log \e_{\omega_{[p+1]}}\prod_{i\leq 2} \e_{z_{[p+1]}^i} \exp \zeta_p
\chi_{i,p+1}\bigl( \omega_{\preceq [p+1]}, z^i_{\preceq [p+1]}\bigr)
\nonumber
\\
\{\mbox{definition (\ref{chij})}\}
= &\
\frac{1}{\zeta_p} \log \e_{\omega_{[p+1]}} \exp \zeta_p
\sum_{i\leq 2}
\chi_{i,p}\bigl( \omega_{\preceq [p+1]}, z^i_{\preceq [p]}\bigr).
\label{Sec5Xp}
\end{align}
For $[p+1]\preceq \alpha$, let us define for $i=1,2$,
\begin{equation}
W_{|\alpha|-1}^i \bigl( \omega_{\preceq [p+1]}, z^i_{\preceq \alpha}\bigr)
=
\exp \zeta_{|\alpha|-1}\Bigl(\chi_{i,|\alpha|}\bigl( \omega_{\preceq [p+1]}, z^i_{\preceq \alpha}\bigr)  - \chi_{i,|\alpha|-1}\bigl( \omega_{\preceq [p+1]}, z^i_{\prec \alpha}\bigr) \Bigr).
\label{Sec5Walphai}
\end{equation}
Comparing this with the definition (\ref{ch31Wp}) and using (\ref{Xichii}) we get that for $[p+1]\prec \alpha$,
\begin{equation}
W_{|\alpha|-1} \bigl( z_{\preceq \alpha} \bigr) 
= 
\prod_{i\leq 2}
W_{|\alpha|-1}^i \bigl( \omega_{\preceq [p+1]}, z^i_{\preceq [j+1]}\bigr).
\label{Sec5Walpha}
\end{equation}
For $\alpha = [p+1]$ this is no longer true, but if we denote
\begin{equation}
\barWp = \barWp \bigl( \omega_{ [p+1]}, z_{\preceq [p]}\bigr) =
\exp \zeta_p\Bigl(\sum_{i\leq 2}\chi_{i,p}\bigl( \omega_{\preceq [p+1]}, z^i_{\preceq [p]}\bigr)  - X_{p}(z_{\preceq [p]}) \Bigr)
\label{Sec5Q}
\end{equation}
then we can write
\begin{equation}
W_p \bigl( z_{\preceq [p+1]} \bigr) 
= 
\barWp \bigl( \omega_{ [p+1]}, z_{\preceq [p]}\bigr)
\prod_{i\leq 2}
W_p^i \bigl( \omega_{\preceq [p+1]}, z^i_{\preceq [p+1]}\bigr).
\label{Sec5Wp}
\end{equation}
If, similarly to (\ref{Sec5Wa}) and (\ref{Sec5WpC}), we denote
\begin{equation}
W_\alpha^i = W_{|\alpha|-1}^i \bigl( z_{\preceq \alpha}\bigr)
\,\,\mbox{ and }\,\,
W_{[p],{C}}^i = \prod_{\alpha\in p([p],{C})} W_\alpha^i
\label{Sec5WpCi}
\end{equation}
then we can rewrite (\ref{Sec5WpC}) as
\begin{equation}
W_{[p],{C}} = \barWp \prod_{i\leq 2} W_{[p],{C}}^{i}.
\label{Sec5Wpr}
\end{equation}
Notice that $\barWp$ does not depend on $z_\beta^i$ for $[p]\prec \beta$, while $W_{[p],{C}}^{1}$ and $\xi_1^{{C}}$ in (\ref{Sec5eq1}) do not depend on $z_\beta^2$ for $[p]\prec \beta$ and $W_{[p],{C}}^{2}$ and $\xi_2^{{C}}$ do not depend on $z_\beta^1$ for $[p]\prec \beta$. This means that if we denote by $\epi$ the expectation in the random variables $z_\beta^i$ for  $[p]\prec \beta$ and denote
\begin{equation}
\eta_i^{{C}} =  \epi\, \xi_i^{{C}} W_{[p],{C}}^{i}
\label{Sec5eqsecond}
\end{equation}
then (\ref{Sec5eq1}) can be rewritten as
\begin{equation}
\e_{\omega_{[p+1]}}  \eta_1^{{C}}  \eta_2^{{C}} \barWp = \e_{\omega_{[p+1]}}   \eta_1^{{C}} \barWp \, \e_{\omega_{[p+1]}}   \eta_2^{{C}} \barWp
\label{Sec5eq2}
\end{equation}
almost surely, because after averaging $\epi$ in $z_\beta^i$ for $[p]\prec \beta$, the only random variable left to be averaged in $\ep$ is $\omega_{[p+1]}$. So far we have just rewritten the equation (\ref{Sec5eq1}) under the induction assumption (\ref{indass}). Now, we will use this to prove the main result of this section.

Let us make the following simple observation: recalling the definition of $A_i^\alpha(\eps)$ in (\ref{Aieps}), if we set all but finite number of random variables $(\pi_i(d))_{d\geq 2}$ to zero, the equation (\ref{Sec5eq2}) still holds almost surely. To see this, first of all, notice that because the random variables $\pi_i(d)$ take any natural value with positive probability, we can set a finite number of them to any values we like in (\ref{Sec5eq2}). For example, for any $D,D'\geq 2$, we can set $\pi_1(d)=\pi_2(d)=n_d$ for $d\leq D$ and set $\pi_1(d)=\pi_2(d)=0$ for $D< d< D'$. The remaining part of the last term in $A_i^\alpha(\eps)$ can be bounded uniformly by
$$
\sum_{d\geq D'} \sum_{k\leq \pi_i(d)} \| \theta_{k,i}^d \|_\infty
\leq
\sum_{d\geq D'} \sum_{k\leq \pi_i(d)}|g_{k,i}^d|,
$$
where, by the assumption (\ref{gsvar}),  we have $\e (g_{k,i}^d)^2 \leq 2^{-d}\epsilon^{\mathrm{pert}},$ which implies that this sum goes to zero almost surely as $D'$ goes to infinity. It follows immediately from this that we can set all but finite number of $\pi_i(d)$ in (\ref{Sec5eq2}) to zero. Moreover, we will set $\pi_1(\lambda K) = \pi_2(\lambda K) = 0$, since the terms coming from the model Hamiltonian will play no role in the proof - all the information we need is encoded into the perturbation Hamiltonian. From now on we will assume that in (\ref{Sec5eq2}), for a given $D\geq 2$,
\begin{equation}
\pi_1(\lambda K) = \pi_2(\lambda K) = 0,\,\,
\pi_1(d)=\pi_2(d)=n_d \,\,\mbox{ for }\,\, d\leq D,\,\, 
\pi_1(d)=\pi_2(d)=0 \,\,\mbox{ for }\,\, d>D.
\label{fixPoisson}
\end{equation}
In addition, let us notice that both sides of (\ref{Sec5eq2}) are continuous functions of the variables $g_i^1$ and $g_{k,i}^d$ for $k\leq n_d$ for $d\leq D$, $i=1,2$. This implies that almost surely over other random variables the equation (\ref{Sec5eq2}) holds for all $g_{k,i}^d$ and, in particular, we can set them to be equal to any prescribed values, 
\begin{equation}
g_1^1 = g_2^1 = g^1,\,\,
g_{k,1}^d=g_{k,2}^d=g_k^d.
\label{fixgs}
\end{equation}
The following is the main result of this section.
\begin{theorem}\label{Sec5Th}
The random variables $\eta_i^{{C}}$ do not depend on $\omega_{[p+1]}$.
\end{theorem}
Here and below, when we say that a function (or random variable) does not depend on a certain coordinate, this means that the function is equal to the average over that coordinate almost surely. In this case, we want to show that
$$
\eta_i^{{C}} = \e_{\omega_{[p+1]}} \eta_i^{{C}}
$$
almost surely.

\medskip
\noindent 
\textbf{Proof of Theorem \ref{Sec5Th}.} Besides the Poisson and Gaussian random variables in (\ref{fixPoisson}), (\ref{fixgs}) and the random variable $\omega_{[p+1]}$ over which we average in (\ref{Sec5eq2}), the random variables $\eta_i^{{C}}$ for $i=1,2$ depend on $\omega_*$, $\omega_{\preceq [p]}$, $(\omega_*^I)_{I\in \II_i}$ and $z_{\preceq [p]}^i$, and $\barWp$ depends on the same random variables for both $i=1,2$. Let us denote
$$
u_i = \bigl( (\omega_*^I)_{I\in \II_i}, z_{\preceq [p]}^i \bigr).
$$
We already stated that for almost all $\omega_*$, $\omega_{\preceq [p]}$, $u_1$ and $u_2$, the equation (\ref{Sec5eq2}) holds for all Poisson and Gaussian random variables fixed as in (\ref{fixPoisson}), (\ref{fixgs}). Therefore, for almost all $\omega_*$, $\omega_{\preceq [p]}$ the equation (\ref{Sec5eq2}) holds for almost all $u_1$, $u_2$ and for all Poisson and Gaussian random variables fixed as in (\ref{fixPoisson}), (\ref{fixgs}). Let us fix any such $\omega_*$, $\omega_{\preceq [p]}$. Then, we can write 
$$
\eta_i^{{C}} = \varphi(u_i,\omega_{[p+1]}) \ \mbox{ and}\
\barWp = \psi(u_1,u_2,\omega_{[p+1]})
$$
for some functions $\varphi$ and $\psi$. These functions depend implicitly on all the random variables we fixed, and the function $\varphi$ is the same for both $\eta_1^{{C}}$ and $\eta_2^{{C}}$ because we fixed all Poisson and Gaussian random variables in (\ref{fixPoisson}), (\ref{fixgs}) to be the same for $i=1,2$. The equation (\ref{Sec5eq2}) can be written as (in the rest of this proof, let us for simplicity of notation write $\omega$ instead of $\omega_{[p+1]}$)
\begin{equation}
\e_{\omega}   \varphi(u_1,\omega) \varphi(u_2,\omega)  \psi(u_1,u_2,\omega)
= 
\prod_{i=1,2}\e_{\omega}   \varphi(u_i,\omega)  \psi(u_1,u_2,\omega)
\label{Sec5eq3}
\end{equation}
for almost all $u_1,u_2$. We want to show that $\varphi(u,\omega)$ does not depend on $\omega$. If we denote
$$
c =: |g_1^1| + |g_2^1| + \sum_{d\leq D} \sum_{k\leq n_d} |g_{k}^d|
$$
then, by (\ref{fixPoisson}) and (\ref{fixgs}), we can bound $A_i^\alpha(\eps)$ in (\ref{Aieps}) by $|A_i^\alpha(\eps)| \leq c$ for $i=1,2$. By induction in (\ref{chij}), $|\chi_{i,|\alpha|}|\leq c$ and, by (\ref{Sec5Xp}), $|X_{p}|\leq 2c$. Therefore, from the definition of $\barWp$ in (\ref{Sec5Q}),
\begin{equation}
e^{-4c}\leq \barWp = \psi(u_1,u_2,\omega) \leq e^{4c}.
\label{Sec5eq5}
\end{equation}
Of course, $|\varphi|\leq 1$. Suppose that for some $\eps>0$, there exists a set $U$ of positive measure such that the variance $\mbox{Var}_{\omega}(\varphi(u,\omega))\geq \eps$ for $u\in U.$ Given $\delta>0$, let $(S_\ell)_{\ell\geq 1}$ be a partition of $L^1([0,1],d\omega)$ such that $\mbox{diam}(S_\ell)\leq \delta$ for all $\ell.$ Let
$$
U_\ell = \bigl\{ u \ |\ \varphi(u,\,\cdot\,) \in S_\ell \bigr\}.
$$
For some $\ell$, the measure of $U\cap U_\ell$ will be positive, so for some $u_1,u_2\in U$,
\begin{equation}
\e_\omega | \varphi(u_1,\omega) - \varphi(u_2,\omega) | \leq \delta.
\label{Sec5eq6}
\end{equation}
The equations (\ref{Sec5eq5}) and (\ref{Sec5eq6}) imply that
$$
\bigl| 
\e_\omega  \varphi(u_1,\omega) \psi(u_1,u_2,\omega)
- \e_\omega  \varphi(u_2,\omega) \psi(u_1,u_2,\omega)
\bigr|
\leq e^{4c}\delta
$$
and, similarly,
$$
\bigl| 
\e_\omega  \varphi(u_1,\omega) \varphi(u_2,\omega)  \psi(u_1,u_2,\omega)
- \e_\omega  \varphi(u_1,\omega)^2 \psi(u_1,u_2,\omega)
\bigr|
\leq e^{4c}\delta.
$$
Since $|\varphi|\leq 1$ and $\e_\omega \psi =1$, the first inequality implies that
$$
\Bigl| 
\prod_{i=1,2}\e_\omega  \varphi(u_i,\omega)  \psi(u_1,u_2,\omega)
- \bigl(\e_\omega  \varphi(u_1,\omega) \psi(u_1,u_2,\omega) \bigr)^2
\Bigr|
\leq
e^{4c}\delta, 
$$
which, together with the second inequality and (\ref{Sec5eq3}), implies
$$
\e_\omega  \varphi(u_1,\omega)^2 \psi(u_1,u_2,\omega)
-
 \bigl(\e_\omega  \varphi(u_1,\omega) \psi(u_1,u_2,\omega) \bigr)^2
 \leq 2 e^{4c}\delta.
$$
The left hand side is a variance with the density $\psi$ and can be written using replicas as
$$
\frac{1}{2} \iint\! \bigl(\varphi(u_1,x)- \varphi(u_1,y)\bigr)^2 \psi(u_1,u_2,x)\psi(u_1,u_2,y)\,dx dy.
$$
By (\ref{Sec5eq5}) and the fact that $u_1\in U$, we can bound this from below by
$$
\frac{1}{2}e^{-8c} \iint\! \bigl(\varphi(u_1,x)- \varphi(u_1,y)\bigr)^2 \,dx dy
=
e^{-8c}\mbox{Var}_{\omega}(\varphi(u_1,\omega)) \geq e^{-8c}\eps.
$$
Comparing lower and upper bounds, $e^{-8c}\eps \leq e^{4c}\delta$, we arrive at contradiction, since $\delta>0$ was arbitrary. Therefore, $\mbox{Var}_{\omega}(\varphi(u,\omega)) = 0$ for almost all $u$ and this finishes the proof.
\qed

\section{A representation formula via properties of the RPC}\label{Sec6label}

Let us summarize what we proved in the previous section. We considered ${C} = \{\alpha_1,\ldots,\alpha_k\}$ for some $[p+1]\preceq \alpha_1,\ldots,\alpha_k \in \Natural^r$ and assumed that $\epi \prod_{\ell\leq k}\sbar_i^{\alpha_\ell}$ does not depend on $\omega_{[p+1]}$. Then, as a consequence of this and the cavity equations, we showed that
\begin{equation}
\eta_i^{{C}} =  \epi\, \xi_i^{{C}} W_{[p],{C}}^{i}
\label{Sec6eqsecond}
\end{equation}
also does not depend on $\omega_{[p+1]}$ almost surely, where 
\begin{equation}
W_{[p],{C}}^i = \prod_{\alpha\in p([p],{C})} W_\alpha^i
\label{Sec6WpCi}
\end{equation}
and $p([p],{C}) = \bigl\{\beta \ \bigr|\ [p+1]\preceq \beta\preceq \alpha, \alpha\in{C} \bigr\}$. Moreover, this holds for Poisson and Gaussian random variables fixed to arbitrary values as in (\ref{fixPoisson}), (\ref{fixgs}). By the assumption (M), the same statement holds in we replace ${C}$ by any subset ${C'}\subseteq C$.

In this section, we will represent the expectation $\epi$ in (\ref{Sec6eqsecond}) with respect to $z_\alpha^i$ for $[p]\prec \alpha$ by using the property of the Ruelle probability cascades in Theorem \ref{Th4label}. Essentially, the expectation in (\ref{Sec6eqsecond}) is of the same type as (\ref{Sec2expectW}) if we think of the vertex $[p]$ as a root. Indeed, we are averaging over random variables indexed by the vertices $[p]\prec \alpha$ which form a tree (if we include the root $[p]$) isomorphic to a tree 
$
\Gamma = \Natural^0\cup \Natural^1 \cup \ldots \cup \Natural^{r-p}
$
of depth $r-p$. We can identify a vertex $[p]\preceq \alpha \in \A$ with the vertex $\{*\} \preceq \gamma \in \Gamma$ such that $\alpha = [p]\gamma$ (for simplicity, we denote by $[p]\gamma$ the concatenation $([p],\gamma)$). Similarly to (\ref{indass}), let us define for $\gamma\in\Natural^{r-p}$,
\begin{equation}
\sbar_I^\gamma = h\bigl( (\omega_\beta)_{\beta\preceq [p]},\omega_{[p+1]}, (\omega_\beta^I)_{\beta\preceq [p]},(\omega_{[p]\beta}^I)_{*\prec \beta \preceq \gamma}  \bigr).
\label{Sec6sialpha}
\end{equation}
Notice a subtle point here: the random variables $\sbar^\gamma_I$ are not exactly the same as $\sbar^\alpha_I$ in (\ref{indass}) for $\alpha = [p]\gamma$. They are exactly the same only if $[p+1]\preceq \alpha$, and in this case we will often write $\sbar_I^\gamma= \sbar_I^{[p]\gamma}.$ Otherwise, if $[p+j]\preceq \alpha$ for $j\geq 2$ then in (\ref{indass}) we plug in the random variable $\omega_{[p+j]}$ instead of $\omega_{[p+1]}$ as we did in (\ref{Sec6sialpha}). The reason for this will become clear soon but, basically, we are going to represent the average $\epi$ in (\ref{Sec6eqsecond}) with respect to $\omega_{[p]\beta}^I$ for $\{*\}\prec \beta$ using the Ruelle probability cascades while $\omega_{[p+1]}$ appears in (\ref{Sec6eqsecond}) on the outside of this average. 

Similarly to (\ref{Aieps}) -- (\ref{Sec3xiialpha}), let us define for $\gamma\in \Natural^{r-p}$,
\begin{align}
A_i^\gamma(\eps) \ =\ & \theta^1(\eps)
+\sum_{2\leq d\leq D} \sum_{k\leq n_d}\theta_{k}^d(\sbar^\gamma_{{1,d,k,i}},\ldots, \sbar^\gamma_{{d-1,d,k,i}},\eps),
\label{Sec6Aieps}
\\
A_i^\gamma \ =\ & \log \Av \exp A_i^\gamma(\eps),
\label{Sec6Aialpha}
\\
\xi_i^\gamma \ =\ & \frac{\Av \eps \exp A_i^\gamma(\eps) }{\exp A_i^\gamma}.
\label{Sec6xiialpha}
\end{align}
The reason why (\ref{Sec6Aieps}) looks different from (\ref{Aieps}) is because we fixed Poisson and Gaussian random variables as in (\ref{fixPoisson}), (\ref{fixgs}), so $\theta^1$ and $\theta_{k}^d$ are defined in terms of $g^1$ and $g_k^d$ in (\ref{fixgs}). Again, let us emphasize one more time that all these definition coincide with the old ones when $[p+1]\preceq [p]\gamma$ or, equivalently, when $[1]\preceq \gamma.$ 

We will keep the dependence on the random variables $(\omega_\beta)_{\beta\preceq [p]},\omega_{[p+1]}, (\omega_\beta^I)_{\beta\preceq [p]}$ implicit and, similarly to (\ref{Sec5xichi}), we will write for $\gamma\in\Natural^{r-p}$,
\begin{equation}
A_i^{\gamma} =\chi_i \bigl((z_{[p]\beta}^i)_{*\prec \beta \preceq \gamma}\bigr). 
\label{Sec6xichi}
\end{equation}
Let $\chi_{i,r}=\chi_i$ and define for $\gamma\in \Gamma\setminus\{*\}$ by decreasing induction on $|\gamma|$,
$$
\chi_{i,p+|\gamma|-1}\bigl( (z_{[p]\beta}^i)_{*\prec \beta \prec \gamma} \bigr)
=
\frac{1}{\zeta_{p+|\gamma|-1}} \log \e_{z_{[p]\gamma}^i} \exp \zeta_{p+|\gamma|-1} \chi_{i,|\alpha|}\bigl( (z_{[p]\beta}^i)_{*\prec \beta \preceq \gamma}\bigr)
$$
and
$$
W_{p+|\gamma|-1}^i \bigl( (z_{[p]\beta}^i)_{*\prec \beta \preceq \gamma} \bigr)
=
\exp \zeta_{p+|\gamma|-1}\Bigl(\chi_{i,p+|\gamma|}\bigl( (z_{[p]\beta}^i)_{*\prec \beta \preceq \gamma}\bigr)  - \chi_{i,p+ |\gamma|-1}\bigl( (z_{[p]\beta}^i)_{*\prec \beta \prec \gamma}\bigr) \Bigr).
$$
For $[1]\preceq \gamma$ these are exactly the same definitions as in (\ref{chij}) and (\ref{Sec5Walphai}), but here we extend these definition to all $\gamma\in \Gamma\setminus \{*\}.$

\smallskip
Let $\tilde{z}_{[p]\beta}^i = (\tilde{\omega}_{[p]\beta}^I)_{I\in\II_i}$ for $\beta \in \Gamma\setminus\{*\}$ be the array generated according to these changes of density along the tree $\Gamma$ as in Section \ref{Sec2ilabel}. Since $[p]$ acts as a root, we do not generate any $\tilde{\omega}_{\beta}^I$ for $|\beta|\leq p$. Similarly to (\ref{Sec4tildas}), we can write for $\gamma\in\Natural^{r-p}$,
\begin{equation}
\xi_i^\gamma = \xi_i \bigl( (z_{[p]\beta}^i)_{*\prec \beta \preceq \gamma}\bigr),\,\,
\tilde{\xi}_i^\gamma = \xi_i \bigl( (\tilde{z}_{[p]\beta}^i)_{*\prec \beta \preceq \gamma} \bigr),
\label{Sec6tildas}
\end{equation}
where we continue to keep the dependence on $(\omega_\beta)_{\beta\preceq [p]},$ $\omega_{[p+1]}$ and $(\omega_\beta^I)_{\beta\preceq [p]}$, as well as Poisson and Gaussian random variables we fixed above, implicit. Given the vertices $[p+1]\preceq \alpha_1,\ldots,\alpha_k \in \Natural^r$ let $[1] \preceq \gamma_1,\ldots,\gamma_k \in \Natural^{r-p}$ be such that $\alpha_\ell = [p]\gamma_\ell.$ Then we can write $\eta_i^{C}$ in (\ref{Sec6eqsecond}) as
\begin{equation}
\eta_i^{C} = \epi\, \xi_i^{\alpha_1}\cdots \xi_i^{\alpha_k} W_{[p],{C}}^{i}
= \e_{*,i}\, \txi_i^{\gamma_1}\cdots \txi_i^{\gamma_k},
\label{Sec6rep1}
\end{equation}
where $\e_{*,i}$ denotes the expectation in $\tilde{z}_{[p]\beta}^i$ for $\beta \in \Gamma\setminus\{*\}$. 

Below we will represent this quantity using the analogue of Theorem \ref{Th4label}. Let $(v_{\gamma})_{\gamma\in\Natural^{r-p}}$ be the weights of the Ruelle probability cascades corresponding to the parameters
\begin{equation}
0<\zeta_{p}<\ldots<\zeta_{r-1}<1,
\end{equation} 
let $(V_\gamma)_{\gamma\in\Natural^{r-p}}$ be their rearrangement as in (\ref{Vs2}) and, similarly to (\ref{tildeVs}), define 
\begin{equation}
\tilde{V}_\gamma
=
\frac{V_\gamma \exp A_i^\gamma}{\sum_{\gamma\in \Natural^{r-p}}V_\gamma \exp A_i^\gamma}.
\label{Sec6tildeVs}
\end{equation}
Theorem \ref{Th4label} can be formulated in this case as follows.
\begin{theorem}\label{Sec6Th4label}
There exists a random bijection $\rho:\Gamma\to\Gamma$ of the vertices of the tree $\Gamma$, which preserves the parent-child relationship, such that
\begin{equation}
\bigl(\tilde{V}_{\rho(\gamma)}\bigr)_{\gamma\in\Natural^{r-p}}
\stackrel{d}{=}
\bigl(V_\gamma\bigr)_{\gamma\in\Natural^{r-p}},\,\,
\bigl(z^i_{[p]\rho(\gamma)} \bigr)_{\gamma\in \Gamma\setminus \{*\}}
\stackrel{d}{=}
\bigl(\tilde{z}^i_{[p]\gamma} \bigr)_{\gamma\in \Gamma\setminus \{*\}}
\label{Sec6Th4eq}
\end{equation}
and these two arrays are independent of each other. 
\end{theorem}
The expectation $\e_{*,i}$ in (\ref{Sec6rep1}) depends on $\gamma_1,\ldots, \gamma_k$ only through their overlaps $(\gamma_\ell\wedge \gamma_{\ell'})_{\ell,\ell'\leq k}$. Above, we made the specific choice $\alpha_\ell = [p]\gamma_\ell$, which implies 
$$
\gamma_\ell\wedge \gamma_{\ell'} = q_{\ell,\ell'} : =\alpha_\ell\wedge \alpha_{\ell'} - p.
$$ 
Now, consider the set of arbitrary configurations with these overlaps,
$$
\CC= \bigl\{ (\gamma_1,\ldots, \gamma_k) \in (\Natural^{r-p})^k \ \bigr|\  \gamma_\ell\wedge \gamma_{\ell'} = q_{\ell,\ell'} \bigr\}.
$$
Let us denote by $\e_*$ the expectation in $(V_\gamma)_{\gamma\in\Natural^{r-p}}$ in addition to $\tilde{z}_{[p]\beta}^i$ for $\beta \in \Gamma\setminus\{*\}$ in the definition of $\e_{*,i}$. We will also denote by $\e_*$ the expectation in $(V_\gamma)_{\gamma\in\Natural^{r-p}}$ and $z_{[p]\beta}^i$ for $\beta \in \Gamma\setminus\{*\}$.  Using Theorem \ref{Sec6Th4label} and arguing as in the proof of Theorem \ref{Sec4Th},
\begin{align}
\e_*
\sum_{(\gamma_1,\ldots, \gamma_k)\in \CC} \tilde{V}_{\gamma_1}\cdots \tilde{V}_{\gamma_k} \,
\xi_i^{\gamma_1}\cdots \xi_i^{\gamma_k}
\ = & \
\e_* \sum_{(\gamma_1,\ldots, \gamma_k)\in \CC} V_{\gamma_1}\cdots V_{\gamma_k} \,
\e_* \txi_i^{\gamma_1}\cdots \txi_i^{\gamma_k}
\nonumber
\\
\ = & \
\eta_i^{C}\
\e_* \sum_{(\gamma_1,\ldots, \gamma_k)\in \CC} V_{\gamma_1}\cdots V_{\gamma_k}.
\label{Sec6repres}
\end{align}
Let us rewrite this equation using a more convenient notation. Let $\sigma^1,\ldots,\sigma^k$ be i.i.d. replicas drawn from $\Natural^{r-p}$ according to the weights $(V_\gamma)_{\gamma\in\Natural^{r-p}}$ and let $\la\,\cdot\, \ra$ denote the average with respect to these weights. If we denote $Q^k = (\sigma^\ell\wedge \sigma^{\ell'})_{\ell,\ell'\leq k}$ and $Q = (q_{\ell,\ell'})_{\ell,\ell'\leq k}$ then we can write
$$
\e_* \sum_{(\gamma_1,\ldots, \gamma_k)\in \CC} V_{\gamma_1}\cdots V_{\gamma_k}
= \e_* \bigl\la \I(Q^k = Q) \bigr\ra = \p(Q^k = Q)
$$
and
$$
\e_*
\sum_{(\gamma_1,\ldots, \gamma_k)\in \CC} \tilde{V}_{\gamma_1}\cdots \tilde{V}_{\gamma_k} \,
\xi_i^{\gamma_1}\cdots \xi_i^{\gamma_k}
=
\e_*
\frac{\bigl\la \xi_i^{\sigma^1}\cdots \xi_i^{\sigma^k} \I(Q^k = Q) \exp \sum_{\ell\leq k} A_i^{\sigma^\ell} \bigr\ra}{\bigl\la \exp A_i^{\sigma} \bigr\ra^k},
$$
and we can rewrite (\ref{Sec6repres}) above as
$$
\e_*
\frac{\bigl\la \prod_{\ell\leq k} \xi_i^{\sigma^\ell} \I(Q^k = Q) \exp \sum_{\ell\leq k} A_i^{\sigma^\ell} \bigr\ra}{\bigl\la \exp A_i^{\sigma} \bigr\ra^k}
=
\eta_i^{C} \p(Q^k = Q).
$$
Notice that this computation also works if we replace each factor $\xi_i^{\gamma_\ell}$ in (\ref{Sec6repres}) by any power $(\xi_i^{\gamma_\ell})^{n_\ell}$ and, in particular, by setting $n_\ell = 0$ or $1$ we get the following. For a subset $S\subseteq \{1,\ldots, k\}$, let us denote $C(S)=\{\gamma_\ell \ | \ \ell\in S\}.$ Then
$$
\e_*
\frac{\bigl\la \prod_{\ell\in S} \xi_i^{\sigma^\ell} \I(Q^k = Q) \exp \sum_{\ell\leq k} A_i^{\sigma^\ell} \bigr\ra}{\bigl\la \exp A_i^{\sigma} \bigr\ra^k}
=
\eta_i^{C(S)} \p(Q^k = Q) .
$$
Furthermore, it will be convenient to rewrite the left hand side using (\ref{Sec6Aialpha}) and (\ref{Sec6xiialpha}) as
$$
\e_*
\frac{\bigl\la \Av \prod_{\ell\in S} \eps_\ell \exp \sum_{\ell\leq k} A_i^{\sigma^\ell}(\eps_\ell)\, \I(Q^k = Q) \bigr\ra}{\bigl\la \Av\, \exp \sum_{\ell\leq k} A_i^{\sigma^\ell}(\eps_\ell) \bigr\ra}.
$$
We showed as a consequence of the assumption (M) that all $\eta_i^{C(S)}$ do not depend on $\omega_{[p+1]}$ and, therefore, we proved the following.
\begin{theorem}\label{Sec6Thend}
Under the assumption (M), for any subset $S\subseteq \{1,\ldots, k\}$,
\begin{equation}
\e_*
\frac{\bigl\la \Av \prod_{\ell\in S} \eps_\ell \exp \sum_{\ell\leq k} A_i^{\sigma^\ell}(\eps_\ell)\, \I(Q^k = Q) \bigr\ra}{\bigl\la \Av\,  \exp \sum_{\ell\leq k} A_i^{\sigma^\ell}(\eps_\ell) \bigr\ra}
\label{Sec6ThendEq}
\end{equation}
does not depend on $\omega_{[p+1]}$ almost surely.
\end{theorem}

\section{Generating Gaussian random fields}\label{Sec7label}

We begin by simplifying (\ref{Sec6Aieps}) further, by taking $n_2=1$ and setting all other $n_d=0$ for $d\geq 3$ except for one, $n_d = n$,
$$
A_i^\gamma(\eps) = \theta^1(\eps)
+ \theta_{1}^2(\sbar^\gamma_{{1,2,1,i}},\eps)
+\sum_{j\leq n}\theta_{j}^d(\sbar^\gamma_{{1,d,j,i}},\ldots, \sbar^\gamma_{{d-1,d,j,i}},\eps).
$$
One can easily see that the definition of $\theta^d$ in (\ref{thetadetx}) satisfies for $\eps\in \{-1,+1\},$
$$
\theta^d(x_1,\ldots,x_{d-1},\eps) =\frac{1+\eps}{2} \log \Bigl(
 1+(e^{g^d}-1) \frac{1+x_{1}}{2}\cdots \frac{1+x_{d-1}}{2} \Bigr)
$$
and, therefore, we can rewrite
$$
A_i^\gamma(\eps)
=
\frac{1+\eps}{2}\Bigl(
g^1 + \log\Bigl(1+(e^{g_1^2}-1)\frac{1+\sbar^\gamma_{1,2,1,i}}{2}\Bigr)
+\sum_{j\leq n} \log\Bigl(1+(e^{g_j^d}-1)\prod_{\ell\leq d-1}\frac{1+\sbar^\gamma_{\ell,d,j,i}}{2}\Bigr)
\Bigr).
$$
At this moment, for simplicity of notation, we will drop some unnecessary indices. We will write $\sbar^\gamma_{i}$ instead of $\sbar^\gamma_{1,2,1,i}$ and, since $d$ is fixed for a moment,
write $\sbar^\gamma_{\ell,j,i}$ instead of $\sbar^\gamma_{\ell,d,j,i}$. Also, we will denote
$$
x_j:= e^{g_j^d}-1 \in (-1,\infty), \,\,
y: = e^{g_1^2}-1 \in (-1,\infty).
$$
Then, we can write
\begin{equation}
A_i^\gamma(\eps)
=
\frac{1+\eps}{2}\Bigl(
g^1 + \log\Bigl(1+y \frac{1+\sbar^\gamma_{i}}{2}\Bigr)
+\sum_{j\leq n} \log\Bigl(1+x_j \prod_{\ell\leq d-1}\frac{1+\sbar^\gamma_{\ell,j,i}}{2}\Bigr)
\Bigr).
\label{SecA1}
\end{equation}
By Theorem \ref{Sec6Thend}, under the assumption (M), the quantities in (\ref{Sec6ThendEq}) do not depend on $\omega_{[p+1]}$ almost surely. In particular, as we discussed above, this almost sure statement can be assumed to hold for all $y,x_j \in (-1,\infty)$ by continuity. We will now take 
\begin{equation}
x_j = x\frac{\eta_j}{\sqrt{n}}
\label{Sec7etas}
\end{equation}
for $x\in (-1,1)$ and independent Rademacher random variables $\eta_j$ and show that, by letting $n\to\infty$, we can replace the last sum in (\ref{SecA1}) by some Gaussian field in the statement of Theorem  \ref{Sec6Thend}. Let us denote
$$
S_{j,i}^\gamma =
\prod_{\ell\leq d-1}\frac{1+\sbar^\gamma_{\ell,j,i}}{2}.
$$
Then with the choice of $x_j = x\eta_j/\sqrt{n}$ we can use Taylor's expansion to write
\begin{equation}
\sum_{j\leq n} \log\Bigl(1+x_j \prod_{\ell\leq d-1}\frac{1+\sbar^\gamma_{\ell,j,i}}{2}\Bigr)
=
\frac{x}{\sqrt{n}} \sum_{j\leq n} \eta_j S_{j,i}^\gamma
- \frac{x^2}{2n} \sum_{j\leq n} (S_{j,i}^\gamma)^2 
+ O(n^{-1/2}).
\label{Sec7napprox}
\end{equation}
The last term $O(n^{-1/2})$ is uniform in all parameters, so it will disappear in (\ref{Sec6ThendEq}) when we let $n$ go to infinity. For the first term, we will use the classical CLT to replace it by Gaussian and for the second term we will use the SLLN, which will produce a term that will cancel out in the numerator and denominator in (\ref{Sec6ThendEq}). However, before we do that, we will need to change the definition of the expectation $\e_*$ slightly. Recall that, by (\ref{Sec6sialpha}),
$$
\sbar^\gamma_{\ell, j,i} 
=
h\bigl( (\omega_\beta)_{\beta \preceq [p]}, \omega_{[p+1]}, (\omega_\beta)_{\beta \preceq [p]}^{\ell, j,i} , (\omega_{[p]\beta}^{\ell, j,i} )_{*\prec \beta \preceq \gamma} \bigr).
$$
In (\ref{Sec6ThendEq}) we already average in the random variables $\omega_{[p]\beta}^{\ell, j,i}$ for $*\prec \beta \preceq \gamma$ but, clearly, the statement of Theorem \ref{Sec6Thend} holds if $\e_*$ also includes the average in $(\omega_\beta^{\ell, j,i})_{\beta\preceq [p]}$ and the Rademacher random variables $\eta_j$ in (\ref{Sec7etas}). From now on we assume this. Note that $\e_*$ still does not  include the average with respect to the random variables $(\omega_\beta^i)_{\beta\preceq [p]}$ that appear in $\sbar_i^\gamma$ in (\ref{SecA1}),
\begin{equation}
\sbar^\gamma_{i} 
=
h\bigl( (\omega_\beta)_{\beta\preceq [p]}, \omega_{[p+1]}, (\omega_\beta^i)_{\beta\preceq [p]}, (\omega_{[p]\beta}^{i} )_{*\prec \beta \preceq \gamma} \bigr).
\label{Sec7sbar}
\end{equation}
Of course, one can not apply the CLT and SLLN in (\ref{Sec6ThendEq}) directly, because there are infinitely many terms indexed by $\gamma\in\Natural^{r-p}$. However, this is not a serious problem because most of the weight of the Ruelle probability cascades $(V_\gamma)$ is concentrated on finitely many indices $\gamma$ and it is not difficult to show that (\ref{Sec6ThendEq}) is well approximated by the analogous quantity where the series over $\gamma$ are truncated at finitely many terms. Moreover, this approximation is uniform over $n$ in (\ref{Sec7napprox}). This is why the representation of $\eta_i^{C(S)}$ in the previous section using the Ruelle probability cascades plays such a crucial role. We will postpone the details until later in this section and first explain what happens in (\ref{Sec7napprox}) for finitely many $\gamma$.

First of all, by the SLLN, for any fixed $(\omega_\beta)_{\beta\preceq [p]}$ and $\omega_{[p+1]}$, 
$$
\lim_{n\to\infty}\frac{1}{n} \sum_{j\leq n} (S_{j,i}^\gamma)^2 
= \e_* (S_{1,i}^\gamma)^2
= \prod_{\ell\leq d-1} \e_* \Bigl(\frac{1+\sbar^\gamma_{\ell,1,i}}{2}\Bigr)^2
= \Bigl(\e_* \Bigl(\frac{1+\sbar^\gamma_{1,1,i}}{2}\Bigr)^2\Bigr)^{d-1}
$$
almost surely. Of course, we can now simplify the notation by replacing $\sbar^\gamma_{1,1,i}$ with $\sbar^\gamma_{i}$ and replacing $\e_*$ by the expectation $\e_i$ with respect to $\omega_\beta^i$ for $\beta\in\A$,
$$
\lim_{n\to\infty}\frac{1}{n} \sum_{j\leq n} (S_{j,i}^\gamma)^2 
= \Bigl(\e_i \Bigl(\frac{1+\sbar^\gamma_{i}}{2}\Bigr)^2\Bigr)^{d-1}
= 2^{2-2d}\Bigl( 1+2\e_i \sbar^\gamma_{i}+\e_i (\sbar^\gamma_{i})^2 \Bigr)^{d-1}.
$$
Lemma \ref{Sec2iLem1} in Section \ref{Sec2ilabel} (for $p=0$ there) yields that $\e_i \sbar^\gamma_{i}$ depends only on $\omega_*$, 
\begin{equation}
q_0(\omega_*): = \e_i \sbar^\gamma_{i},
\label{Sec7q0}
\end{equation}
and, since $\e_i (\sbar^\gamma_{i})^2$ clearly does not depend on $\gamma$, we can take $[1]\preceq \gamma$, in which case 
$$
\e_i (\sbar^\gamma_{i})^2 = \e_i (\sbar^{[p]\gamma}_{i})^2 = R_{[p]\gamma,[p]\gamma} = q_r.
$$
This means that for any $\gamma\in \Natural^{r-p}$,
\begin{equation}
\lim_{n\to\infty}\frac{1}{n} \sum_{j\leq n} (S_{j,i}^\gamma)^2 = 2^{2-2d}\bigl(1+2q_0(\omega_*)+q_r \bigr)^{d-1}
\label{Sec7LLN}
\end{equation}
almost surely. So, in the limit, these terms will cancel out in (\ref{Sec6ThendEq}) -- at least when we truncate the summation over $\gamma$ to finitely many $\gamma$, as we shall do below.

\smallskip
Next, let us look at the first sum in (\ref{Sec7napprox}) for $\gamma\in F$ for a finite set $F\subset \Natural^{r-p}$. By the classical multivariate CLT (applied for a fixed $(\omega_\beta)_{\beta\preceq [p]}$ and $\omega_{[p+1]}$),
\begin{equation}
\frac{1}{\sqrt{n}} \sum_{j\leq n} \eta_k \bigl(S_{j,i}^\gamma\bigr)_{\gamma\in F}
\stackrel{d}{\longrightarrow}
(g^\gamma)_{\gamma\in F},
\label{Sec7CLT}
\end{equation}
where $(g^\gamma)_{\gamma\in F}$ is a centered Gaussian random vector with the covariance 
$$
\e g^\gamma g^{\gamma'}
=
\e_* S_{1,i}^\gamma S_{1,i}^{\gamma'}
=
\Bigl(
\e_i \frac{1+\sbar_i^\gamma}{2}\cdot\frac{1+\sbar_i^{\gamma'}}{2}
\Bigr)^{d-1}
=
\Bigl( \frac{1}{4}\bigl(1+\e_i \sbar^\gamma_{i}+\e_i \sbar^{\gamma'}_{i}+\e_i \sbar^\gamma_{i} \sbar^{\gamma'}_{i} \bigr)\Bigr)^{d-1}.
$$
First of all, as above $\e_i \sbar^\gamma_{i} = \e_i \sbar^{\gamma'}_{i} = q_0(\omega_*).$ To compute $\e_i \sbar^\gamma_{i} \sbar^{\gamma'}_{i}$, we need to consider two cases. First, suppose that $\gamma\wedge\gamma'\geq 1.$ Since $\e_i \sbar^\gamma_{i} \sbar^{\gamma'}_{i}$ clearly depends only on $\gamma\wedge\gamma'$, we can suppose that $[1]\preceq \gamma,\gamma'$, in which case $\sbar^{\gamma}_{i} = \sbar^{[p]\gamma}_{i}$, $\sbar^{\gamma'}_{i} = \sbar^{[p]\gamma'}_{i}$ and
$$
\e_i \sbar^\gamma_{i} \sbar^{\gamma'}_{i}
=
\e_i \sbar^{[p]\gamma}_{i} \sbar^{[p]\gamma'}_{i}
= R_{[p]\gamma,[p]\gamma'} = q_{p+\gamma\wedge\gamma'}.
$$
In the second case, $\gamma\wedge\gamma' = 0$, when computing $\e_i \sbar^\gamma_{i} \sbar^{\gamma'}_{i}$ we can first average $\sbar_i^\gamma$ with respect to $ (\omega_{[p]\beta}^{i} )_{*\prec \beta \preceq \gamma}$ and $\sbar_i^{\gamma'}$ with respect to $ (\omega_{[p]\beta}^{i} )_{*\prec \beta \preceq \gamma'}$, since these are independent. However, by Lemma \ref{Sec2iLem1}, both of these averages do not depend on $\omega_{[p+1]}$. This means that, after taking these averages,  we can replace $\omega_{[p+1]}$ in (\ref{Sec7sbar}) by $\omega_{[p+j]}$ if $[j]\preceq \gamma$, and the same for $\gamma'$. As a result, we can again write
$$
\e_i \sbar^\gamma_{i} \sbar^{\gamma'}_{i}
=
\e_i \sbar^{[p]\gamma}_{i} \sbar^{[p]\gamma'}_{i}
= R_{[p]\gamma,[p]\gamma'} = q_p = q_{p+\gamma\wedge\gamma'}
$$
almost surely. If we introduce the notation
$$
c_{j}(\omega_*) = \frac{1+2q_0(\omega_*) + q_{p+j}}{4}
\,\,\mbox{ for }\,\,
j=0,\ldots, r-p
$$
then we proved that the covariance is given by
\begin{equation}
\e g^\gamma g^{\gamma'}
=
\e_* S_{1,i}^\gamma S_{1,i}^{\gamma'}
=
c_{\gamma\wedge\gamma'}(\omega_*)^{d-1}.
\label{Sec7cov}
\end{equation}
Let us show right away that 
\begin{equation}
0\leq c_0(\omega_*) < \ldots < c_{r-p}(\omega_*).
\label{Sec7cs}
\end{equation}
In particular, this means that the covariance in (\ref{Sec7cov}) is increasing with $\gamma\wedge\gamma'$ and the Gaussian field $(g^\gamma)_{\gamma\in\Natural^{r-p}}$ is the familiar field that accompanies the Ruelle probability cascades in the pure $d$-spin non-diluted models. Of course, the only statement that requires a proof is the following.
\begin{lemma}
The inequality $c_0(\omega_*)\geq 0$ holds almost surely.
\end{lemma} 
\textbf{Proof.} First of all, let us notice that, by Lemma \ref{Sec2iLem1}, we can also write the definition of $q_0(\omega_*)$ in (\ref{Sec7q0}) as
$
q_0(\omega_*) = \e_i \sbar^\alpha_{i}
$
for any vertex $\alpha\in\Natural^r$, where $\e_i$ denotes the expectation in the random variables $\omega_\beta^i$ for $\beta\in\A.$ Let $\e_i'$ denote the expectation with respect to $\omega_\beta^i$ for $\beta\in\A\setminus \{*\}$, excluding the average in $\omega_*^i$. Again, by Lemma \ref{Sec2iLem1}, we can write $\e_i' \sbar^\alpha_{i}$ as $q_0(\omega_*,\omega_*^i)$. If we take any $\alpha,\alpha'\in \Natural^{r}$ such that $\alpha\wedge\alpha'=0$ then
$$
q_0(\omega_*)^2 = (\e_i \e_i'\sbar^\alpha_{i})^2
\leq  \e_i (\e_i'\sbar^\alpha_{i})^2 = \e_i (\e_i'\sbar^\alpha_{i} \e_i'\sbar^{\alpha'}_{i})
=
\e_i  \sbar^\alpha_{i} \sbar^{\alpha'}_{i}
=q_{\alpha\wedge\alpha'} = q_0
$$
almost surely. Therefore
$$
4c_0(\omega_*)\geq 1-2\sqrt{q_0} +q_0 = (1-\sqrt{q_0})^2\geq 0.
$$
This finishes the proof.
\qed

\medskip
From now on let $\e_*$ also include the expectation in the Gaussian field $(g^\gamma)_{\gamma\in\Natural^{r-p}}$ (conditionally on $\omega_*$) with the covariance (\ref{Sec7cov}). Let us denote by
\begin{equation}
a_i^\gamma(\eps)
=
\frac{1+\eps}{2}\Bigl(
g^1 + \log\Bigl(1+y \frac{1+\sbar^\gamma_{i}}{2}\Bigr)
+x g^\gamma
\Bigr)
\label{Sec7a1}
\end{equation}
the quantity that will replace $A_i^\gamma(\eps)$ in (\ref{SecA1}) in the limit $n\to\infty$. We are ready to prove the following.
\begin{theorem}\label{Sec7Th}
Under the assumption (M), for any subset $S\subseteq \{1,\ldots, k\}$,
\begin{equation}
\e_*
\frac{\bigl\la \Av \prod_{\ell\in S} \eps_\ell \exp \sum_{\ell\leq k} a_i^{\sigma^\ell}(\eps_\ell)\, \I(Q^k = Q) \bigr\ra}{\bigl\la \Av\, \exp \sum_{\ell\leq k} a_i^{\sigma^\ell}(\eps_\ell) \bigr\ra}
\label{Sec7ThEq}
\end{equation}
does not depend on $\omega_{[p+1]}$ almost surely.
\end{theorem}
\textbf{Proof.} We only need to show that the quantity in (\ref{Sec6ThendEq}) with $A_i^\gamma(\eps)$ as in (\ref{SecA1}) with the choice of $x_j$ as in (\ref{Sec7etas}) converges to (\ref{Sec7ThEq}), where $\e_*$ was redefined above (\ref{Sec7sbar}) to include the average over $(\omega_\beta)_{\beta\preceq [p]}^{\ell, j,i}$ and Rademacher random variables $\eta_j$ in (\ref{Sec7etas}), as well as the Gaussian field $(g^\gamma)_{\gamma\in\Natural^{r-p}}$ with the covariance (\ref{Sec7cov}).

Let $(V_m')_{m\geq 1}$ be the RPC weights $(V_\gamma)_{\gamma\in\Natural^{r-p}}$ arranged in the decreasing order. For some fixed large $M\geq 1$, let us separate the averages $\la\, \cdot\, \ra$ in the numerator and denominator in (\ref{Sec6ThendEq}) into two sums over $V_m'$ for $m\leq M$ and for $m>M$ (for each replica). Let us denote the corresponding averages by $\la\, \cdot\, \ra_{\leq M}$ and $\la\, \cdot\, \ra_{>M}$. Let
\begin{align*}
a \ =&\ 
\bigl\la \Av \prod_{\ell\in S} \eps_\ell \exp \sum_{\ell\leq k} A_i^{\sigma^\ell}(\eps_\ell)\, \I(Q^k = Q) \bigr\ra_{\leq M},
\\
b \ =&\ 
\bigl\la \Av \prod_{\ell\in S} \eps_\ell \exp \sum_{\ell\leq k} A_i^{\sigma^\ell}(\eps_\ell)\, \I(Q^k = Q) \bigr\ra_{> M},
\\
c \ =&\ 
\bigl\la \Av\,  \exp \sum_{\ell\leq k} A_i^{\sigma^\ell}(\eps_\ell) \bigr\ra_{\leq M},
\\
d \ =&\ 
\bigl\la \Av\,  \exp \sum_{\ell\leq k} A_i^{\sigma^\ell}(\eps_\ell) \bigr\ra_{> M}.
\end{align*}
Note that $|a|\leq c$, $|b|\leq d$ and
$$
\Bigl|\frac{a+b}{c+d} - \frac{a}{c}\Bigr|
=
\Bigl|\frac{bc-ad}{c(c+d)}\Bigr|
\leq
\Bigl|\frac{b}{c+d}\Bigr|+\Bigl|\frac{d}{c+d}\Bigr|
\leq
\frac{2d}{c+d}.
$$
This means that the difference between (\ref{Sec6ThendEq}) and
\begin{equation}
\e_* \frac{a}{c}
=
\e_*
\frac{\bigl\la \Av \prod_{\ell\in S} \eps_\ell \exp \sum_{\ell\leq k} A_i^{\sigma^\ell}(\eps_\ell)\, \I(Q^k = Q) \bigr\ra_{\leq M}}{\bigl\la \Av\,  \exp \sum_{\ell\leq k} A_i^{\sigma^\ell}(\eps_\ell) \bigr\ra_{\leq M}}
\label{Sec7trunc}
\end{equation}
can be bounded by $2\e_* d/(c+d)$. By (\ref{Sec7napprox}), the SLLN in (\ref{Sec7LLN}) and CLT in (\ref{Sec7CLT}), with probability one, (\ref{Sec7trunc}) converges to
\begin{equation}
\e_*
\frac{\bigl\la \Av \prod_{\ell\in S} \eps_\ell \exp \sum_{\ell\leq k} a_i^{\sigma^\ell}(\eps_\ell)\, \I(Q^k = Q) \bigr\ra_{\leq M}}{\bigl\la \Av\,  \exp \sum_{\ell\leq k} a_i^{\sigma^\ell}(\eps_\ell) \bigr\ra_{\leq M}}
\label{Sec7limM}
\end{equation}
as $n\to\infty$. Next, we show that $\e_* d/(c+d)$ is small for large $M$, uniformly over $n$. Since $d/(c+d)\in [0,1]$, it is enough to show that $d$ is small and $c$ is not too small with high probability. To show that $d$ is small, we will use Chebyshev's inequality and show that $\e_* d$ is small. By Jensen's inequality,
$$
d  = 
\bigl\la \Av\,  \exp \sum_{\ell\leq k} A_i^{\sigma^\ell}(\eps_\ell) \bigr\ra_{> M}
=
\bigl\la \Av\,  \exp A_i^{\sigma}(\eps) \bigr\ra_{> M}^k
\leq
\bigl\la \Av\,  \exp k\, A_i^{\sigma}(\eps) \bigr\ra_{> M}.
$$
If we denote $\delta_M = \e \sum_{m>M} V_m'$ then, since the weights $(V_\gamma)$ and the random variables in $A_i^\gamma(\eps)$ are independent,
$$
\e_*d \leq  
\delta_M
\sup_{\gamma}  \Av\, \e_* \exp k\, A_i^{\gamma}(\eps).
$$
Using that $\log(1+t)\leq t$, we can bound $A_i^{\gamma}(\eps)$ with the choice (\ref{Sec7etas}) by
$$
A_i^{\gamma}(\eps)
\leq |g^1| + \log\bigl(1+|y|\bigr) + \frac{1+\eps}{2}\frac{x}{\sqrt{n}} \sum_{j\leq n} \eta_j S_{j,i}^\gamma.
$$
Using that, for a Rademacher random variable $\eta$, we have $\e e^{\eta t} = \ch(t) \leq e^{t^2/2}$ we get that
$$
\e_\eta \exp  k\frac{1+\eps}{2}\frac{x}{\sqrt{n}} \sum_{j\leq n} \eta_j S_{j,i}^\gamma 
\leq e^{k^2/2}
$$
and, therefore,
$$
\sup_{\gamma}  \Av\, \e_* \exp k\, A_i^{\gamma}(\eps)
\leq
c_k:=\exp\Bigl( k|g^1| + k\log\bigl(1+|y|\bigr) +\frac{k^2}{2} \Bigr).
$$
We showed that $\e_* d \leq c_k \delta_M$, and this bound does not depend on $n$. Since $\delta_M$ is small for $M$ large, $d$ is small with high probability uniformly over $n.$ On the other hand, to show that $c$ is not too small, we simply bound it from below by
$$
c\geq \bigl(V_\gamma\, \Av \exp A_i^{\gamma}(\eps)\bigr)^k
$$
for $\gamma$ corresponding to the largest weight, $V_1' = V_\gamma.$ The weight $V_1'$ is strictly positive and its distribution does not depend on $n$. Also, using (\ref{Sec7napprox}), we can bound $A_i^{\gamma}(\eps)$ from below by
$$
A_i^{\gamma}(\eps)
\geq 
 -|g^1| + \log\bigl(1-|y|\bigr) - \Bigl|\frac{1}{\sqrt{n}} \sum_{j\leq n} \eta_j S_{j,i}^\gamma \Bigr| - L 
$$
for some absolute constant $L$. Even though the index $\gamma$ here is random, because it corresponds to the largest weight $V_1'$, we can control this quantity using Hoeffding's inequality for Rademacher random variables conditionally on other random variables to get
$$
\p\Bigl(
\Bigl|\frac{1}{\sqrt{n}} \sum_{j\leq n} \eta_j S_{j,i}^\gamma \Bigr| \geq t
\Bigr)
\leq 2e^{-t^2/2}.
$$
Therefore, for any $\delta>0$ there exists $\Delta>0$ (that depends on $|g^1|$, $|y|$, $k$ and the distribution of $V_1'$) such that $\p(c\geq \Delta)\geq 1-\delta.$ All together, we showed that $\e_* d/(c+d)$ is small for large $M$, uniformly over $n$.

To finish the proof, we need to show that (\ref{Sec7limM}) approximates (\ref{Sec7ThEq}) for large $M$. Clearly, this can be done by the same argument (only easier) that we used above to show that (\ref{Sec7trunc}) approximates (\ref{Sec6ThendEq}).
\qed

\section{Final consequences of the cavity equations}

Theorem \ref{Sec7Th} implies that, under the assumption (M),
\begin{equation}
\e_*
\frac{\bigl\la \Av \prod_{\ell\leq k}(1+ \eps_\ell) \exp \sum_{\ell\leq k} a_i^{\sigma^\ell}(\eps_\ell)\, \I(Q^k = Q) \bigr\ra}{\bigl\la \Av\, \exp \sum_{\ell\leq k} a_i^{\sigma^\ell}(\eps_\ell) \bigr\ra}
\label{Sec8Eq1}
\end{equation}
does not depend on $\omega_{[p+1]}$ almost surely. This follows from (\ref{Sec7ThEq}) by multiplying out
$\prod_{\ell\leq k}(1+ \eps_\ell).$
Using that for $\eps\in\{-1,+1\}$,
$$
\exp t \frac{1+\eps}{2} = 1+ (e^t-1)\frac{1+\eps}{2}, 
$$
one can see from (\ref{Sec7a1}) that
$$
2\Av \exp
a_i^\gamma(\eps)
=
1+e^{g^1+xg^\gamma}
\Bigl(1+y \frac{1+\sbar^\gamma_{i}}{2}\Bigr)
$$
and\smallskip
$$
\Av (1+\eps) \exp a_i^\gamma(\eps)
=
e^{g^1+xg^\gamma}
\Bigl(1+y \frac{1+\sbar^\gamma_{i}}{2}\Bigr).
\medskip
$$
If for simplicity the notation we denote $z:=e^{g^1}\in(0,\infty)$ and
\begin{equation}
S_i^\gamma = \frac{1+\sbar^\gamma_{i}}{2}
\label{Sec8Sgamma}
\end{equation}
then (\ref{Sec8Eq1}) can be written, up to a factor $2^k z^k$, as
$$
\e_*
\frac{\bigl\la \prod_{\ell\leq k} \exp(xg^{\sigma^\ell})(1+yS_i^{\sigma^\ell})\, \I(Q^k = Q) \bigr\ra}{\bigl\la 1+z \exp(xg^{\sigma})(1+yS_i^{\sigma}) \bigr\ra^k}.
$$
As before, the statement that this quantity does not depend on $\omega_{[p+1]}$ almost surely holds for all $x\in (-1,1)$, $y>-1$ and $z>0$, by continuity. Therefore, if we take the derivative with respect to $z$ and then let $z\downarrow 0$ then the quantity we get (up to a factor $-k$),
$$
\e_* \Bigl\la \prod_{\ell\leq k+1} \exp(xg^{\sigma^\ell})(1+yS_i^{\sigma^\ell})\, \I(Q^k = Q) \Bigr\ra,
$$
does not depend on $\omega_{[p+1]}$ almost surely. This is a polynomial in $y$ of order $k+1$ and if we take the derivative $\partial^{k+1}/\partial y^{k+1}$ we get that
$$
\e_* \Bigl\la \prod_{\ell\leq k+1} \exp(x g^{\sigma^\ell}) \,S_i^{\sigma^\ell} \, \I(Q^k = Q) \Bigr\ra,
$$
does not depend on $\omega_{[p+1]}$ almost surely. Let us now take the expectation $\e_g$ with respect to the Gaussian field $(g^\gamma)$. By (\ref{Sec7cov}), 
$$
\e_g \exp x\sum_{\ell\leq k+1} g^{\sigma^\ell}
=
\exp \frac{x^2}{2}\sum_{\ell,\ell'\leq k+1} c_{\sigma^\ell\wedge \sigma^{\ell'}}(\omega_*)^{d-1}
$$
and, therefore, the quantity
$$
\e_* \Bigl\la \prod_{\ell\leq k+1} S_i^{\sigma^\ell} \, \I(Q^k = Q) 
\exp \frac{x^2}{2}\sum_{\ell,\ell'\leq k+1} c_{\sigma^\ell\wedge \sigma^{\ell'}}(\omega_*)^{d-1}
\Bigr\ra
$$
does not depend on $\omega_{[p+1]}$ almost surely. Notice that because of the indicator $\I(Q^k = Q)$, all the overlaps $\sigma^\ell\wedge\sigma^{\ell'}$ are fixed for $\ell,\ell' \leq k$, so the factor
$$
\exp \frac{x^2}{2}\sum_{\ell,\ell'\leq k} c_{\sigma^\ell\wedge \sigma^{\ell'}}(\omega_*)^{d-1}
$$
can be taken outside of $\e_* \la\,\cdot \,\ra$ and cancelled out, yielding that
$$
\e_* \Bigl\la \prod_{\ell\leq k+1} S_i^{\sigma^\ell} \, \I(Q^k = Q) 
\exp \frac{x^2}{2}\sum_{\ell\leq k} c_{\sigma^\ell\wedge \sigma^{k+1}}(\omega_*)^{d-1}
\Bigr\ra
$$
does not depend on $\omega_{[p+1]}$ almost surely. Taking the derivative with respect to $x^2/2$ at zero gives that
\begin{equation}
\e_* \Bigl\la \prod_{\ell\leq k+1} S_i^{\sigma^\ell} \, \I(Q^k = Q) 
\sum_{\ell \leq k} c_{\sigma^\ell\wedge \sigma^{k+1}}(\omega_*)^{d-1}
\Bigr\ra
\label{Sec9polyd}
\end{equation}
does not depend on $\omega_{[p+1]}$ almost surely. We proved this statement for a fixed but arbitrary $d\geq 3$, but it also holds for $d=1$ by setting $x=0$ in the previous equation. Let us take arbitrary $f_j$ for $j=0,\ldots,r-p$ and consider a continuous function $f$ on $[0,1]$ such that $f(c_j(\omega_*)) = f_j.$ Approximating this function by polynomials, (\ref{Sec9polyd}) implies that
$$
\e_* \Bigl\la \prod_{\ell\leq k+1} S_i^{\sigma^\ell} \, \I(Q^k = Q) 
\sum_{\ell \leq k} \sum_{j=0}^{r-p}f_j \I(\sigma^\ell\wedge \sigma^{k+1}=j)
\Bigr\ra
$$
does not depend on $\omega_{[p+1]}$ for all $(f_j)$ almost surely. Taking the derivative in $f_j$ shows that
$$
\e_* \Bigl\la \prod_{\ell\leq k+1} S_i^{\sigma^\ell} \, \I(Q^k = Q) 
\sum_{\ell \leq k} \I(\sigma^\ell\wedge \sigma^{k+1}=j)
\Bigr\ra
$$
does not depend on $\omega_{[p+1]}$ almost surely and, therefore,
$$
\e_* \Bigl\la \prod_{\ell\leq k+1} S_i^{\sigma^\ell} \, \I(Q^k = Q) 
\sum_{\ell \leq k} \I(\sigma^\ell\wedge \sigma^{k+1}\geq j)
\Bigr\ra
$$
does not depend on $\omega_{[p+1]}$ almost surely for all $j=0,\ldots, r-p$. Let us now express this quantity as a linear combination over all possible overlap configurations that the new replica $\sigma^{k+1}$ can form with the old replicas $\sigma^1,\ldots,\sigma^k$.

Given a $k\times k$ overlap constraint matrix $Q=(q_{\ell,\ell'})_{\ell,\ell'\leq k}$, let $\EE(Q)$ be the set of admissible extensions $Q'=(q_{\ell,\ell'}')_{\ell,\ell'\leq k+1}$ of $Q$ to $(k+1)\times(k+1)$ overlap constraint matrices. In other words, $q_{\ell,\ell'}' = q_{\ell,\ell'}$ for $\ell,\ell'\leq k$, and there exists $\gamma_1,\ldots,\gamma_{k+1}\in \Natural^{r-p}$ such that $\gamma_\ell\wedge \gamma_{\ell'} = q_{\ell,\ell'}'$ for $\ell,\ell'\leq k+1.$ If we denote
\begin{equation}
n_j(Q') = \sum_{\ell \leq k} \I(q_{\ell,k+1}'\geq j)
\label{Sec8nj}
\end{equation}
and denote $Q^{k+1} = (\sigma^\ell\wedge\sigma^{\ell'})_{\ell,\ell'\leq k+1}$ then we showed that
$$
\sum_{Q'\in\EE(Q)} n_j(Q')\,
\e_* \Bigl\la \prod_{\ell\leq k+1} S_i^{\sigma^\ell} \, \I(Q^{k+1} = Q') 
\Bigr\ra
$$
does not depend on $\omega_{[p+1]}$ almost surely for all $j=0,\ldots, r-p$. Since the RPC weights $(V_\gamma)$ are independent of $(\sbar_i^\gamma)$, we can rewrite this as
\begin{equation}
\sum_{Q'\in\EE(Q)} n_j(Q') \p(Q^{k+1} = Q')\,  M(Q')
\label{Sec8at}
\end{equation}
does not depend on $\omega_{[p+1]}$ almost surely for all $j=0,\ldots, r-p$, where 
\begin{equation}
M(Q') = \e_* \prod_{\ell\leq k+1} S_i^{\gamma_\ell}
\end{equation}
for any $\gamma_1,\ldots,\gamma_{k+1}\in \Natural^{r-p}$ such that $\gamma_\ell\wedge \gamma_{\ell'} = q_{\ell,\ell'}'$ for $\ell,\ell'\leq k+1.$

Recall that this statement was proved under the induction assumption (M) in Section \ref{Sec2ilabel}, so let us express (\ref{Sec8at}) in the notation of Section \ref{Sec2ilabel} and connect everything back to the assumption (M), which we repeat one more time. Given $[p+1]\preceq \alpha_1,\ldots,\alpha_k \in \Natural^r$, we assumed that:
\begin{enumerate}
\item[(M)] for any subset $S\subseteq \{1,\ldots, k\}$, $\epi \prod_{\ell\in S}\sbar_i^{\alpha_\ell}$ does not depend on $\omega_{[p+1]}$.
\end{enumerate} 
If similarly to (\ref{Sec8Sgamma}) we denote, for $\alpha\in\Natural^r$,
\begin{equation}
S_i^\alpha = \frac{1+\sbar^\alpha_{i}}{2}
\label{Sec8Salpha}
\end{equation}
then the assumption (M) is, obviously, equivalent to
\begin{enumerate}
\item[(M)] for any subset $S\subseteq \{1,\ldots, k\}$, $\epi \prod_{\ell\in S}S_i^{\alpha_\ell}$ does not depend on $\omega_{[p+1]}$.
\end{enumerate} 
Let us define $[1] \preceq \gamma_1,\ldots,\gamma_{k}\in \Natural^{r-p}$ by $\alpha_\ell = [p]\gamma_\ell$, and let $Q$ be the overlap matrix
\begin{equation}
q_{\ell,\ell'} = \gamma^\ell\wedge \gamma^{\ell'}
\label{Sec8Qk}
\end{equation}
The assumption (M) depends on $\alpha_1,\ldots,\alpha_k$ only through this matrix $Q$, so one should really view it as a statement about such $Q$. Fix $1\leq j \leq r-p$ and consider any $Q'\in\EE(Q)$ such that $n_j(Q')\not = 0.$ Since
$$
n_j(Q') \not = 0 \Longleftrightarrow \max_{\ell\leq k} q_{\ell,k+1}'\geq j,
$$
this means that $q_{\ell,k+1}' \geq j\geq 1$ for some $\ell\leq k$. Since our choice of $\alpha_1,\ldots,\alpha_k \in \Natural^r$ was such that $q_{\ell,\ell'} = \gamma_\ell\wedge \gamma_{\ell'}\geq 1$ for all $\ell,\ell'\leq k$, this also implies that $q_{\ell,\ell'}'\geq 1$ for all $\ell,\ell'\leq k+1$. In particular, we can find $[1] \preceq \gamma_{k+1}\in \Natural^{r-p}$ such that $\gamma_\ell\wedge \gamma_{k+1} = q_{\ell,k+1}'$ for $\ell\leq k$. Let $\alpha_{k+1} = [p]\gamma_{k+1}.$ Recall that, whenever $\alpha = [p]\gamma$ for $[1] \preceq \gamma\in \Natural^{r-p}$, the definitions (\ref{Sec6sialpha}) and (\ref{indass}) imply that $\sbar^\alpha_i = \sbar_i^\gamma$. Therefore, in this case, we can rewrite the definition of $M(Q')$ below (\ref{Sec8at}) as
\begin{equation}
M(Q') = \e_* \prod_{\ell\leq k+1} S_i^{\gamma_\ell}
= \epi \prod_{\ell\leq k+1} S_i^{\alpha_\ell}.
\end{equation}
Let us summarize what we proved.
\begin{theorem}\label{Sec8Th}
If the matrix $Q$ defined in (\ref{Sec8Qk}) satisfies the assumption (M) then
\begin{equation}
\sum_{Q'\in\EE(Q)} n_j(Q') \p(Q^{k+1} = Q')\,  M(Q')
\label{Sec8ThEq}
\end{equation}
does not depend on $\omega_{[p+1]}$ almost surely for all $j=1,\ldots, r-p$.
\end{theorem}

\section{Main induction argument}\label{Sec9label}

Finally, we will now use Theorem \ref{Sec8Th} to prove our main goal, Theorem \ref{Sec6iTh1} in Section \ref{Sec2ilabel}. To emphasize that our inductive proof will have a monotonicity property (M), we can rephrase Theorem \ref{Sec6iTh1} as follows.
\begin{theorem}\label{Sec9Th1}
Under the assumption (\ref{indass}), for any $k\geq 1$, any $[p+1]\preceq \alpha_1,\ldots,\alpha_k \in \Natural^r$ and any $S\subseteq \{1,\ldots, k\}$, the expectation $\epi \prod_{\ell\in S} S_i^{\alpha_\ell}$ with respect to $(\omega_\beta^i)_{[p+1]\preceq \beta}$ does not depend on $\omega_{[p+1]}$ almost surely.
\end{theorem}

\begin{figure}[t]
 \centering
 \psfrag{a_0}{ {\small root $[p+1]$ at depth zero}}\psfrag{a_2}{ {\small level $(r-p-1)$ of $\alpha\in \Natural^r$}}\psfrag{AcodeR}{ {\small level $(r-p)$ of multiplicities}}
 \includegraphics[width=0.45\textwidth]{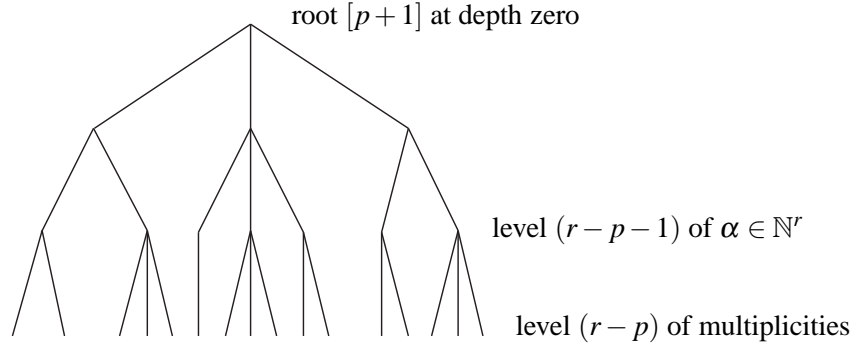}
 \caption{\label{Fig1} Representing a configuration $\alpha_1,\ldots,\alpha_k$ by a tree.}
 \end{figure}

It is much easier to describe the proof if we represent a configuration $[p+1]\preceq \alpha_1,\ldots,\alpha_k \in \Natural^r$ not by a matrix $Q = (\gamma_\ell\wedge\gamma_{\ell'})$ with $\alpha_\ell = [p]\gamma_\ell$ but by a subtree of $\A$ growing out of the vertex $[p+1]$ with branches leading to the leaves $\alpha_1,\ldots,\alpha_k$, and with an additional layer encoding their multiplicities (see Fig. \ref{Fig1}). If we think of $[p+1]$ as a root of this subtree being at depth zero, then the leaves are at depth $r-p-1.$ However, the multiplicity of any particular vertex $\alpha$ in the set $\{\alpha_1,\ldots,\alpha_k\}$ can be greater than one, so we will attach that number of children to each vertex $\alpha$ to represent multiplicities, so the depth of the tree will be $r-p$. Whenever we say that we remove a leaf $\alpha$ from the tree, we mean that we remove one multiplicity of $\alpha$. Notice that removing a leaf from the tree also removes the path to that leaf, of course, keeping the shared paths leading to other leaves that are still there.

We will say that this \textit{tree is good} if 
$$
M(Q) = \epi \prod_{\ell\leq k} S_i^{\alpha_\ell}
$$ 
does not depend on $\omega_{[p+1]}$.

\begin{figure}[t]
 \centering
 \psfrag{a_0}{ {\small root $[p+1]$ at depth zero}}\psfrag{a_1}{\small $n_j$ fixed children}\psfrag{a_2}{ {\small designated vertex at depth $j$}}\psfrag{AcodeR}{ {\small new path leading to some $\alpha\in \Natural^r$}}\psfrag{T_i}{\small ${\cal T}$} 
 \includegraphics[width=0.4\textwidth]{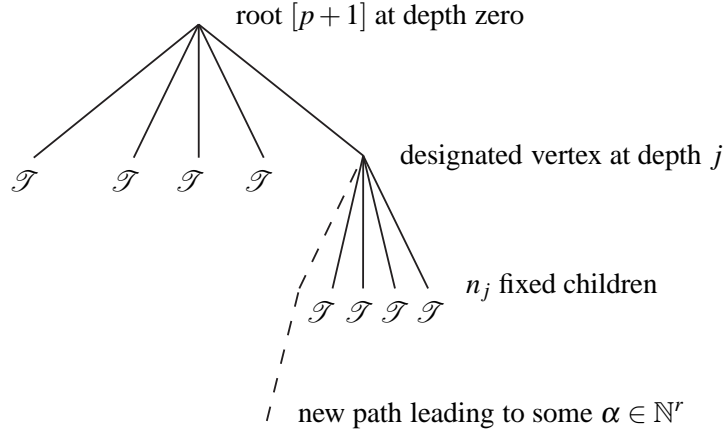}
 \caption{\label{Fig2} Illustrating property ${\cal N}_j$. Solid lines represent the subtree ${\cal T}_j$ and $n_j$ children of the designated vertex at depth $j$. Each $\cal T$ represents an arbitrary tree leading to some set of leaves with their multiplicities. Dashed line represents a new path from the designated vertex to a new vertex $\alpha\in \Natural^r$, which has multiplicity one.}
\end{figure}

First we are going to prove the following property, illustrated in Fig. \ref{Fig2}, that we will call property ${\cal N}_j$ for $j=0,\ldots, r-p-1.$ Let us consider an arbitrary configuration of paths leading from the root $[p+1]$ to some set of vertices at depth $j$. Let us call this part of the tree ${\cal T}_j$, which is now fixed. We pick one designated vertex at depth $j$ (right-most vertex at depth $j$ in Fig. \ref{Fig2}). To all other vertices at depth $j$ we attach arbitrary trees $\cal T$ leading to some arbitrary finite sets of leaves in $\Natural^r$ and their multiplicities. We will use the same generic notation $\cal T$ to represent an arbitrary tree, even though they can all be different. The designated vertex has some fixed number of children, say $n_j$, and to each of these children we also attach an arbitrary tree $\cal T$. Property ${\cal N}_j$ will be the following statement.
\begin{enumerate}
\item[(${\cal N}_j$)] Fix any ${\cal T}_j$ and the number of children $n_j$ of a designated vertex. Suppose that all trees that we just described are good (this means for all choices of  trees $\cal T$, possibly empty). Then any tree obtained by adding a single new path leading from a designated vertex at depth $j$ to some new vertex $\alpha\in \Natural^r$ with multiplicity one (as in Fig. \ref{Fig2}) is also good.
\end{enumerate}

\begin{lemma}
For any $j=0,\ldots, r-p-1,$ any ${\cal T}_j$ and $n_j$, the property ${\cal N}_j$ holds.
\end{lemma}
\textbf{Proof.} Let us fix any particular choice of trees $\cal T$ attached to non-designated vertices at depth $j$ and $n_j$ children of the designated vertex. By the assumption in property ${\cal N}_j$, this tree, as well as any tree obtained  by removing a finite number of leaves, is good. This precisely means that the assumption (M) holds for this tree, or for the sets of leaves with their multiplicities encoded by this tree. Let $Q$ be the matrix described above Theorem \ref{Sec8Th} corresponding to this set of leaves. Then, Theorem \ref{Sec8Th} implies that
\begin{equation}
\sum_{Q'\in\EE(Q)} n_j(Q') \p(Q^{k+1} = Q')\,  M(Q')
\label{Sec9Lem1eq}
\end{equation}
does not depend on $\omega_{[p+1]}$. Obviously, each $Q'\in \EE(Q)$ corresponds to a new tree constructed by adding one more new vertex $\alpha$ to our tree or increasing the multiplicity of some old vertex by one. Moreover, if $n_j(Q') \not = 0$ then the overlap $\alpha \wedge \alpha_\ell$ with one of the old vertices $\alpha_\ell$ should be greater or equal than $j$. This means that this new vertex will be attached to the tree somewhere at depth $j$ or below. One of the possibilities is described in Fig. \ref{Fig2}, when $\alpha$ is attached by a new path to the designated vertex at depth $j$. All other possibilities -- attaching $\alpha$ below one of the non-designated vertices at depth $j$ or below one of the $n_j$ children of the designated vertex -- would simply modify one of the trees $\cal T$ in Fig. \ref{Fig2}. But such a modification results in a good tree, by the assumption in property ${\cal N}_j$. Since the sum in (\ref{Sec9Lem1eq}) is a linear combination of all these possibilities, the term corresponding to adding a new path as in Fig. \ref{Fig2} must be good, which finishes the proof. 
\qed

\begin{figure}[t]
 \centering
 \psfrag{a_0}{ {\small root $[p+1]$ at depth zero}}\psfrag{a_1}{\small single path to $\alpha\in \Natural^r$}\psfrag{a_2}{ {\small designated vertex at depth $j$}}\psfrag{T_i}{\small ${\cal T}$} 
 \includegraphics[width=0.35\textwidth]{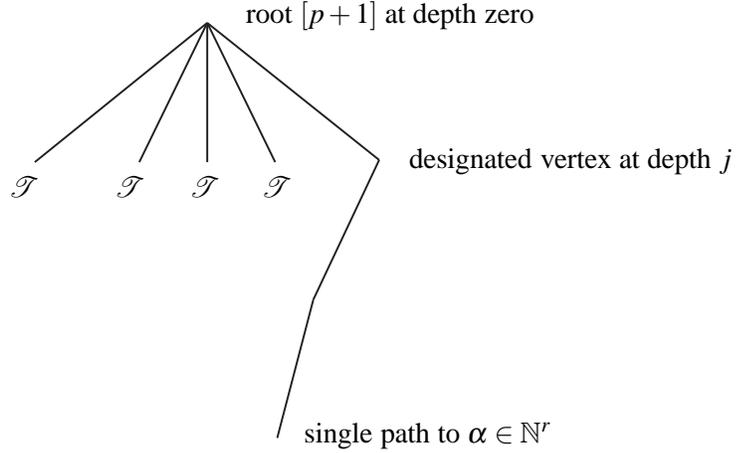}
 \caption{\label{Fig3} Illustrating property ${\cal P}_j$. Solid lines represent the subtree ${\cal T}_j$ and one path from a designated vertex at depth $j$ to a leaf $\alpha\in \Natural^r$ with multiplicity one. Property ${\cal P}_j$ allows to replace this single path by an arbitrary tree $\cal T$.}
\end{figure}

\medskip
\noindent
Next we will prove another property that we will denote $\PP_j$ for $j=0,\ldots,r-p-1$, described in Fig. \ref{Fig3}. As in Fig. \ref{Fig2}, we consider an arbitrary configuration ${\cal T}_j$ of paths leading from the root $[p+1]$ to some set of vertices at depth $j$ and we pick one designated vertex among them. To all other vertices at depth $j$ we attach arbitrary trees $\cal T$, while to the designated vertex we attach a single path leading to some vertex $\alpha\in\Natural^r$ with multiplicity one. Property $\PP_j$ will be the following statement. 
\begin{enumerate}
\item[($\PP_j$)] Suppose that the following holds for any fixed tree ${\cal T}_j$ up to depth $j$. Suppose that any tree as in Fig. \ref{Fig3} is good, as well as any tree obtained by removing any finite number of leaves from this tree. Then any tree obtained by replacing a single path below the designated vertex at depth $j$ by an arbitrary tree $\cal T$ is also good.
\end{enumerate}
We will now prove the following.
\begin{lemma}
Property $\PP_j$ holds for any $j=0,\ldots,r-p-1$.
\end{lemma}
\textbf{Proof.} First of all, notice that property $\PP_{r-p-1}$ follows immediately from property ${\cal N}_{r-p-1}$. In property ${\cal N}_{r-p-1}$, arbitrary trees $\cal T$ below non-designated vertices at depth $r-p-1$ represent their arbitrary multiplicities, the trees $\cal T$ below the children of the designated vertex are empty, and the multiplicity of the designated vertex is $n_j$. Property ${\cal N}_{r-p-1}$ then implies that we can increase this multiplicity by one to $n_j+1.$ Starting from multiplicity one and using this repeatedly, we can make this multiplicity arbitrary. This is exactly the property $\PP_{r-p-1}$.

\begin{figure}[t]
 \centering
 \psfrag{a_0}{ {\small root $[p+1]$ at depth zero}}\psfrag{a_1}{$\ldots$}\psfrag{a_2}{ {\small designated vertex $\delta_j$ at depth $j$}}\psfrag{T_i}{\small ${\cal T}$} 
 \includegraphics[width=0.4\textwidth]{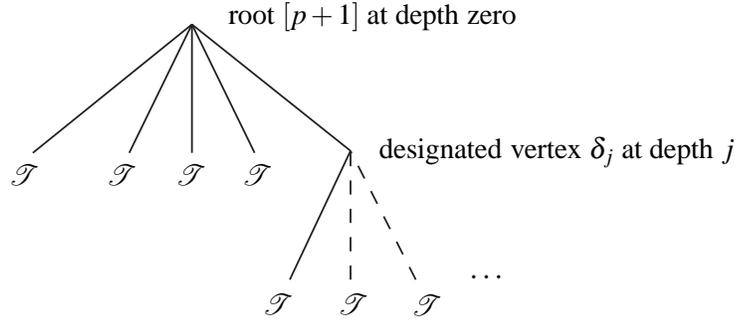}
 \caption{\label{Fig4} Illustrating proof of property ${\cal P}_j$. First we replace single path  in Fig. \ref{Fig3} by arbitrary tree below the child of the designated vertex using ${\cal P}_{j+1}$. Then we iteratively add a new path using property ${\cal N}_j$ and then replace this path below depth $j+1$ be arbitrary tree using ${\cal P}_{j+1}$.}
\end{figure}

Next, we are going to show that property $\PP_{j+1}$ implies property $\PP_j$. The proof of this is illustrated in Fig. \ref{Fig4}. Given any tree as in Fig. \ref{Fig3}, let us denote by $\delta_j$ the designated vertex at depth $j$. Consider the subtree ${\cal T}_{j+1}$ up to depth $j+1.$ It forms the same pattern, with the child of $\delta_j$ playing a role of the designated vertex at depth $j+1.$ Therefore, by property $\PP_{j+1}$ we can replace the single path below this vertex by an arbitrary tree $\cal T$. By property ${\cal N}_j$, if we attach another path to $\delta_j$, the resulting new tree is good. Then we can again treat the child of $\delta_j$ along this new path as a designated vertex at depth $j+1$, apply property $\PP_{j+1}$ and replace the path below this vertex by an arbitrary tree. If we continue to repeatedly use property ${\cal N}_j$ to attach another path to $\delta_j$ and then use property $\PP_{j+1}$ to replace the part of this path below depth $j+1$ by an arbitrary tree, we can create an arbitrary tree below $\delta_j$, and this tree is good by construction. This is precisely property $\PP_j$, so the proof is completed by decreasing induction on $j.$
\qed

\medskip
\noindent
Finally, this implies Theorem \ref{Sec9Th1} (and Theorem \ref{Sec6iTh1}). As we explained in Section \ref{Sec2ilabel}, by induction on $p$ this implies Theorem \ref{Th2} .

\smallskip
\noindent
\textbf{Proof of Theorem \ref{Sec9Th1}.} By Lemma \ref{Sec2iLem1}, the tree consisting of one path from $[p+1]$ (at depth zero) to some vertex $\alpha\in\Natural^{r}$ (at depth $r-p-1$) with multiplicity one is good. Using property $\PP_0$ implies that arbitrary finite tree is good, which finishes the proof.
\qed

\begin{thebibliography}{999}

\bibitem{AC} Aizenman, M., Contucci, P.: On the stability of the quenched state in mean-field spin-glass models. J. Statist. Phys. \textbf{92}, no. 5-6, 765--783 (1998)

\bibitem{AS2} Aizenman, M., Sims, R., Starr, S.L.: An extended variational principle for the SK spin-glass model. Phys. Rev. B. \textbf{68}, 214403 (2003)

\bibitem{Aldous} Aldous, D.: Representations for partially exchangeable arrays of random variables.  J. Multivariate Anal.  {11}, no. 4, 581--598 (1981) 

\bibitem{AA} Arguin, L.-P., Aizenman, M.: On the structure of quasi-stationary competing particles systems. Ann. Probab. \textbf{37}, no. 3, 1080--1113 (2009)

\bibitem{Austin} Austin, T.: Exchangeable random measures. To appear, Ann. Inst. Henri Poincar\'e Probab. Stat., arXiv:1302.2116 (2013)

\bibitem{AP} Austin, T., Panchenko, D.: A hierarchical version of the de Finetti and Aldous-Hoover representations. To appear in Probab. Theory Relat. Fields, arXiv:1301.1259 (2013)

\bibitem{Bolthausen} Bolthausen, E., Sznitman, A.-S.: On Ruelle's probability cascades and an abstract cavity method. Comm. Math. Phys. \textbf{197}, no. 2, 247--276 (1998) 

\bibitem{FL} Franz, S., Leone, M.: Replica bounds for optimization problems and diluted spin systems. J. Statist. Phys. 111, no. 3-4, 535--564 (2003)

\bibitem{GG} Ghirlanda, S., Guerra, F.: General properties of overlap probability distributions in disordered spin systems. Towards Parisi ultrametricity.  J. Phys. A  {31}, no. 46, 9149--9155 (1998) 

\bibitem{Guerra} Guerra, F.: Broken replica symmetry bounds in the mean field spin glass model. Comm. Math. Phys. {\bf 233}, no. 1, 1--12 (2003)

\bibitem{Hoover2} Hoover, D. N.: Row-column exchangeability and a generalized model for probability. Exchangeability in probability and statistics (Rome, 1981), pp. 281--291, North-Holland, Amsterdam-New York (1982) 

\bibitem{Kallenberg} Kallenberg, O.: On the representation theorem for exchangeable arrays. {J. Multivariate Anal.}, {30}, no. 1, 137--154 (1989)

\bibitem{Mezard} M\'ezard, M., Parisi, G.: The Bethe lattice spin glass revisited. Eur. Phys. J. B Condens. Matter Phys. {20}, no. 2, 217--233 (2001)

\bibitem{PT} Panchenko, D., Talagrand, M.: Bounds for diluted mean-fields spin glass models. Probab. Theory Related Fields 130,  no. 3, 319--336 (2004) 

\bibitem{PUltra} Panchenko, D.: The Parisi ultrametricity conjecture. Ann. of Math. (2) {177}, no. 1, 383--393 (2013)

\bibitem{PPF} Panchenko, D.: The Parisi formula for mixed $p$-spin models. To appear in Ann. of Probab., arXiv:1112.4409 (2011)

\bibitem{SKmodel} Panchenko, D.: The Sherrington-Kirkpatrick Model. Springer Monographs in Mathematics. Springer-Verlag, New York (2013)

\bibitem{Pspins} Panchenko, D.: Spin glass models from the point of view of spin distributions.  Ann. of Probab.  {41}, no. 3A, 1315--1361 (2013)

\bibitem{HEPS} Panchenko, D.: Hierarchical exchangeability of pure states in mean field spin glass models. Preprint, arXiv:1307.2207 (2013)

\bibitem{1RSB} Panchenko, D.: Structure of $1$-RSB asymptotic Gibbs measures in the diluted $p$-spin models. J. Statist. Phys. 155, no. 1, 1--22 (2014)

\bibitem{Parisi79} Parisi, G.: Infinite number of order parameters for spin-glasses. Phys. Rev. Lett. \textbf{43}, 1754--1756 (1979)

\bibitem{Parisi} Parisi, G.: A sequence of approximate solutions to the S-K model for spin glasses. J. Phys. A \textbf{13}, L-115 (1980) 

\bibitem{Ruelle} Ruelle, D.: A mathematical reformulation of Derrida's REM and GREM.  Comm. Math. Phys.  {\bf 108},  no. 2, 225--239 (1987)

\bibitem{SK} Sherrington, D., Kirkpatrick, S.: Solvable model of a spin glass. Phys. Rev. Lett. {\bf 35}, 1792--1796 (1975)

\bibitem{TPF} Talagrand, M.: The Parisi formula. Ann. of Math. (2) \textbf{163}, no. 1, 221--263 (2006)

\bibitem{SG} Talagrand, M.: Mean-Field Models for Spin Glasses. Ergebnisse der Mathematik und ihrer Grenzgebiete. 3. Folge A Series of Modern Surveys in Mathematics, Vol. 54, 55. Springer-Verlag (2011)



\end{thebibliography}
\end{document}